\documentclass[a4paper,12pt]{article}         
\usepackage{amsmath,amsfonts,amsthm,amssymb}  
\usepackage[utf8]{inputenc} 
\usepackage[T1]{fontenc}    
									   
\DeclareMathOperator*{\argmin}{arg\,min} 

\usepackage{graphicx,caption}            
\usepackage{color}                       
\usepackage{hyperref}                    
\usepackage{grffile}

\usepackage{makeidx}                  
\usepackage{acronym}
\usepackage{mathrsfs}
\usepackage{float}
\usepackage{vruler}                    
\usepackage{palatino, url, multicol}
\usepackage{setspace}                  
\usepackage{txfonts}

\usepackage{verbatim}

\usepackage{tikz}
\usepackage{tikz-3dplot}
\usetikzlibrary{shapes.geometric, arrows, calc, 3d}

\usepackage[numbers,sort&compress]{natbib}   
\usepackage{import}

\usepackage{caption}

\usepackage{listings}

\usepackage{comment}

\restylefloat{figure}

\title{Flexible Foil Mesh Generation for Spatial \\ 
Focal-Body Modeling of a Spherical Mirror}

\author{Netzer Moriya}
\date{}

\begin{document}

\maketitle

\begin{abstract}
We present a novel application of the Flexible Foil Mesh Generation (FFMG) method to model 
the $3D$ Focal Body generated by a spherical mirror collecting light from an infinitely distant source 
on its optical axis. 
The study addresses the challenge of accurately representing highly concave structures formed by the 
focusing effect. Through theoretical analysis and numerical simulations, we demonstrate the 
effectiveness of the FFMG method in capturing the intricate geometry of the Focal Body, 
with implications for computational geometry, $3D$ reconstruction, and optical system modeling. 
\end{abstract}

\section{Introduction}
\subsection{Background}
Accurate $3D$ shape modeling is needed in fields such as computational geometry, $3D$ reconstruction, 
and spatial data analysis. 
These fields often require the ability to represent both convex and concave structures within a point 
cloud, particularly in applications where the geometric intricacies of an object are critical. 
One such example is in optical systems, where various components including shaped mirrors,
are used to focus light rays onto a confined volume. 
The resulting shape, referred to here as a Focal Body (FB), may prove highly concave and 
poses significant challenges in accurate geometric modeling.
One option to capture the geometric structure of such highly concentrated intensity regions can be by
referring to close proximities of rays as intersecting points that form a spatially confined point 
cloud~\cite{2024Vaara, klimes1994network}.
Traditional methods in computational geometry, such as Convex Hulls or Alpha Shapes, are effective 
in approximating the outer boundary of a point cloud but often prove insufficient when comes to accurately 
capturing the complex concave features inherent in Focal Bodies. 
These limitations can lead to inaccuracies in applications involving optical design, where precise modeling 
of light paths and focal regions is crucial.
The Flexible Foil Mesh Generation (FFMG) method, originally developed to handle highly concave 
structures within point clouds~\cite{Moriya:2025FlexableFoil, Moriya:2024formconvexhullconcavity}, 
offers a promising solution to these challenges. By tending to minimize the 
surface area while ensuring all (or partial) points within the cloud are enclosed, the FFMG method 
provides a more accurate representation of complex geometries, making it particularly suited for applications 
involving concave shapes and ensuring closed surfaces. This study aims to explore the application of FFMG in modeling the $3D$ 
Focal Body formed by a spherical mirror, demonstrating its effectiveness in capturing the 
intricate concave structures that traditional methods struggle to represent.

\subsection{Objective}

The primary objective of this study is to apply the FFMG method
for the precise modeling of the three-dimensional Focal Body formed by the focusing of collimated light 
along the optical axis of a spherical mirror. 
The principal aim is to assess the capability of FFMG in accurately representing the complex and highly 
concave geometry characteristic of the Focal Body.
In this study, we define the \emph{Focal Body} as the three-dimensional spatial region formed by the intersection 
envelope of reflected rays from a spherical mirror, given an infinitely distant collimated light source. 
This definition relies purely on geometrical optics and does not explicitly incorporate intensity variations or 
wavefront interference effects. The FB is constructed by tracking the spatial density of ray intersections within a 
predefined volume, resulting in a highly concave surface representation. This approach ensures a computationally efficient 
and analytically tractable model, particularly for applications in optical design and computational geometry.

\vspace{0.4cm}

Specifically, this paper seeks to:
\begin{enumerate}
    \item Develop a detailed methodology for applying the FFMG method to the problem of modeling a Focal Body.
    \item Validate the method through numerical simulations.
    \item Analyze the performance of FFMG in terms of accuracy, computational efficiency, and its ability to 
	handle concave structures.
\end{enumerate}

\subsection{Related Work}

$3D$ surface reconstruction is a fundamental problem in computational geometry that involves constructing 
a surface that approximates the shape of an object based on a set of sample points, 
typically obtained from real-world data. This process is crucial in various applications, 
including computer graphics, medical imaging, and reverse engineering. 
The challenge lies in accurately capturing both the global shape and fine details, particularly 
in the presence of noise, outliers, and varying point densities.

Historically, several approaches have been developed to tackle $3D$ surface reconstruction. 
One of the earliest and most widely used methods is the \emph{Convex Hull}~\cite{Loffler:2010}, 
which constructs the 
smallest convex shape that encloses all points in the dataset. While simple and computationally 
efficient, the Convex Hull fails to capture concave regions, making it inadequate for more complex shapes.

To address this limitation, more sophisticated techniques have been developed, 
including \emph{Delaunay Triangulation}~\cite{su1996delaunay} 
and \emph{Voronoi Diagrams}~\cite{aurenhammer1991voronoi}, which are often used to generate a 
mesh that represents the surface. These methods offer a balance between computational efficiency 
and accuracy, but they can struggle with highly concave shapes or non-uniform point distributions.

Another notable method is \emph{Alpha Shapes}~\cite{edelsbrunner1994alpha}, which extends the concept 
of the Convex Hull by 
introducing a parameter that controls the level of detail captured in the shape. This allows 
for the representation of concavities, but the choice of the alpha parameter can be challenging 
and may require tuning based on the specific dataset. Also, the method may not guarantee a final 
closed surface.

Optical systems often involve complex geometries, including highly concave surfaces, where traditional 
modeling techniques may fail to meet capturing the precise shape and focusing properties.

One of the most fundamental optical components is the \emph{spherical mirror}, which focuses parallel 
light rays originated from an infinitely distant source at along the optical axis, onto a confined 
spatial volume (i.e., non-single point) in the vicinity of its optical axis. 
The shape of the focused light region, referred here as the \emph{Focal Body}, is a three-dimensional 
structure that reflects the concave shape of the mirror's surface. 
Accurate modeling of such focal bodies for more complex mirror geometries is essential for predicting 
and optimizing mirror performance, particularly in applications such as telescopes, laser systems, 
and imaging devices.

Traditional methods for geometric modeling in optical systems, including ray tracing and wavefront 
analysis, provide insights into light propagation and focusing behavior. However, these methods often 
require precise surface models to accurately predict how light interacts with the optical elements. 
When it comes to representing the Focal Body, traditional surface reconstruction techniques, 
may not fully capture the concave geometry, leading to potential inaccuracies in optical simulations.

A more detailed discussion on essential methods for generating $3D$ mesh surfaces from point clouds can 
be found in \cite{Moriya:2025FlexableFoil,Moriya:2024TheLargest}.

\section{Theoretical Background}

\subsection{A Brief Analysis of Light Distribution Reflected by a Spherical Mirror}

\subsubsection{Problem Setup}
Consider a three-dimensional spherical mirror with an infinitely distant light source positioned along 
the primary optical axis of the mirror. Given the infinite distance of the source, the incoming rays 
are parallel to each other and to the primary optical axis. It is assumed that when parallel rays 
strike a spherical mirror, they reflect according to the law of reflection, where the angle of 
incidence equals the angle of reflection. 

Rays that strike the center of the mirror (on the optical axis) will reflect directly back along the same 
path. However, rays that strike off-center will reflect at an angle. For simplicity, we assume the light 
is monochromatic and disregard its wave nature.

\subsubsection{Ray Behavior Near the Focal Point}

The focal point of the mirror is located at half the radius of curvature. This is a consequence of the mirror equation:

\begin{equation}
\begin{aligned}
\frac{1}{f} &= \frac{1}{d_o} + \frac{1}{d_i} \
&= \lim_{d_o \to \infty} \left(\frac{1}{d_o} + \frac{1}{d_i}\right) \
&= \frac{1}{d_i} \
\therefore f &= d_i
\end{aligned}
\end{equation}

where $d_o$ is the object distance (considered infinite for parallel rays) and $d_i$ is the image distance. \\
From the geometry of spherical mirrors, we know that for parallel rays, the image forms at a 
distance of  $\frac{R}{2}$ from the mirror's surface. 
Therefore, $d_i = \frac{R}{2}$, and consequently, 

\begin{equation} \label{eq:mirror_equation}
f = \frac{R}{2}
\end{equation}

Paraxial rays, which are close to the optical axis, converge near the focal point. For these rays, the 
small-angle approximation $\sin{\theta} \approx \theta$ (where $\theta$ is the angle of incidence or 
reflection) holds true, resulting in reflection through or very near to the focal point, which lies at a 
distance $f$ from the mirror's vertex.

In contrast, rays farther from the optical axis intersect the axis at points closer to the mirror due to 
spherical aberration. Spherical aberration occurs because non-paraxial rays, which make larger angles 
with the optical axis, do not adhere closely to the small-angle approximation. 
The intersection point $z$ on the optical axis for a ray striking the mirror at a distance $y$ from 
the optical axis can be formulated by:

\begin{equation}
z = \frac{R \left[1 - 2\left(\frac{y^2}{R^2}\right)\right]}{\sqrt{1 - \frac{y^2}{R^2}}}
\end{equation}

where \( R \) is the radius of curvature, and \( y \) is the height of the incident ray above the optical axis. To simplify this expression, we expand the square root term using a Taylor series for small \( y/R \):

\begin{equation}
\sqrt{1 - \frac{y^2}{R^2}} \approx 1 - \frac{y^2}{2R^2} - \frac{y^4}{8R^4}
\end{equation}

Substituting this into the expression for \( z \) and simplifying, we get:

\begin{equation}
x \approx R \left[1 - \frac{y^2}{R^2} + \frac{y^4}{4R^4}\right]
\end{equation}

By substituting \eqref{eq:mirror_equation}, we obtain the higher-order approximation:

\begin{equation}
z \approx f - \frac{y^4}{16f^3}  \label{eq:TylorExp}
\end{equation}

This approximation shows that non-paraxial rays intersect the optical axis closer to the mirror 
than the paraxial focal point, with the deviation increasing rapidly with \( y^4 \), reflecting the 
impact of spherical aberration.

While the approximation in Equation \eqref{eq:TylorExp} accurately describes small-angle deviations within the paraxial region, its 
validity diminishes for rays striking the mirror at large angles (i.e., marginal rays with high $\sin{\theta} \approx \theta$, no longer holds, 
and higher-order terms in the expansion become non-negligible). Additionally, the assumption that all rays follow a perfect geometric path 
without considering diffraction or wavefront curvature may introduce further discrepancies in real-world optical applications.

To improve precision in modeling large-angle behavior, alternative numerical strategies such as higher-order aberration corrections, 
or wavefront-based approaches could be employed. These methods would account for higher-order deviations that are otherwise neglected in 
the current approximation. 
Future work could explore such refinements, particularly in applications requiring high-precision optical simulations where marginal ray 
deviations contribute significantly to focal distortions.

To understand the impact of spherical aberration, we may consider the longitudinal spherical aberration, 
which represents the difference in the focal lengths between paraxial rays (close to the optical axis) 
and non-paraxial rays (farther from the axis). 
In spherical mirrors, rays near the optical axis converge at the focal length $f$, 
where \( R \) is the radius of curvature. However, non-paraxial rays intersect the optical axis at points 
closer to the mirror due to spherical aberration.

For a spherical mirror with radius of curvature \( R \), the focal length for paraxial rays is $f$. 
Non-paraxial rays, striking the mirror at a height \( y \) from the optical axis, experience 
spherical aberration, shifting the focal point.

The image distance \( d_i \) for such rays can be approximated using:

\begin{equation}
d_i = \frac{R \cos 2\theta}{\cos \theta}
\end{equation}

where \( \theta \approx \frac{y}{R} \). Expanding \( \cos \theta \) and \( \cos 2\theta \) for small \( y/R \):

\begin{equation}
\cos \theta \approx 1 - \frac{y^2}{2R^2}, \quad \cos 2\theta \approx 1 - 2\frac{y^2}{R^2}
\end{equation}

leads to:

\begin{equation}
d_i \approx R \left(1 - \frac{y^2}{R^2}\right)
\end{equation}

The longitudinal spherical aberration \( \Delta d \) is:

\begin{equation}
\Delta d = d_i - f = -\frac{y^2}{2R}
\end{equation}

Thus, the aberration, growing quadratically with \( y \), is (see also below in \eqref{eq:y2/R2}):

\begin{equation}
\text{Longitudinal aberration} \approx \frac{y^2}{2R}  \label{eq:y2/R2a}
\end{equation}

This approximation shows that the aberration grows quadratically with the distance \( y \) from the 
optical axis, indicating that off-axis rays are increasingly affected. 
The quadratic relationship means that as \( y \) increases, the error in the convergence point 
grows significantly, leading to a spread in the focal point and a corresponding reduction in 
image sharpness. Minimizing this aberration is crucial in optical design to ensure high-quality 
image formation.

The paraxial approximation is central to many optical designs because it ensures that rays converge at a 
common focal point, producing sharp images. However, as rays deviate from the paraxial region, 
aberrations like spherical aberration become significant, affecting the quality of the image.

For a more detailed mathematical treatment of these concepts, including the derivation of the expressions 
for spherical aberration and a deeper understanding of ray tracing through spherical mirrors, refer 
to~\cite{fowles1975, born1999}.

\subsubsection{Formation of the Caustic Surface}

The envelope of reflected rays forms a caustic surface, which is a region where light intensity is 
significantly higher due to the concentration of rays. In the context of a spherical mirror, 
this caustic surface is a result of spherical aberration, where rays that are not paraxial (i.e., 
rays far from the optical axis) do not converge precisely at the focal point but instead form a 
curved surface with a cusp.

For a spherical mirror, the caustic surface has a cusp at the focal point and extends from $\frac{3}{4}R$ 
to $R$ from the mirror's vertex, where $R$ denotes the radius of curvature\footnote{$\frac{3}{4}R$ is approximately
the location where the marginal rays (rays hitting the edge of the mirror) cross the optical axis. However, the caustic 
extends all the way to R}. 
The shape of this caustic can be understood more deeply by considering the reflection of rays at different distances from the optical 
axis.

The reflected rays' envelope can be mathematically described using the concept of geometric optics. 
For a given incident ray at a height $y$ from the optical axis, the reflected ray's angle relative to 
the optical axis varies with $y$, leading to a variation in the intersection points along the optical 
axis. The cusp at the focal point is where the reflected rays are most tightly concentrated, and as 
we move away from this point, the rays spread out, forming the extended part of the caustic.

The curvature of the caustic surface near the focal point can be analyzed using catastrophe 
theory~\cite{arnold1992catastrophe}, 
which provides a mathematical framework for understanding the cusp and fold structures in such optical 
systems. Specifically, the caustic is an example of a "cusp catastrophe," where the intensity of light 
diverges near the cusp.

The distance $z$ along the optical axis where the caustic 
forms can be expressed as a function of the ray height $y$~\cite{berry1980}:

\begin{equation}
z(y) = f - \frac{y^4}{16f^3}
\end{equation}

This relationship indicates that as $y$ increases, the intersection point $z$ moves closer to the mirror, 
contributing to the formation of the caustic surface. The distance from the mirror's vertex to the point 
on the optical axis where the caustic terminates corresponds to the mirror's radius of curvature, $R$.

\subsubsection{Spherical Aberration Effects}

Rays originating from the outer portions of the mirror, also known as marginal rays, focus closer to the 
mirror than paraxial rays, leading to the formation of a blur circle instead of a perfect point focus. 
This phenomenon is known as spherical aberration and occurs because the curvature of a spherical 
mirror does not perfectly focus all incoming parallel rays to a single point.

\begin{enumerate}

   \item \textbf{Quantitative Analysis of Spherical Aberration:}
   The difference in focal points between paraxial rays (which are close to the optical axis) and marginal 
   rays (which are farther from the optical axis) is the primary cause of spherical aberration. For a 
   spherical mirror, the focal length for rays at a distance \( y \) from the optical axis can be 
   approximated as (see also above in \eqref{eq:y2/R2a}):

   \begin{equation}
   f(y) = f_0 \left(1 - \frac{y^2}{2R^2}\right) \label{eq:y2/R2}
   \end{equation}

   where:
   \begin{itemize}
       \item \( f_0 = \frac{R}{2} \) is the focal length for paraxial rays,
       \item \( R \) is the radius of curvature of the mirror,
       \item \( y \) is the distance from the optical axis.
   \end{itemize}

   This expression shows that as \( y \) increases, the effective focal length \( f(y) \) decreases, 
   causing marginal rays to focus closer to the mirror than paraxial rays.

	\item \textbf{Formation of the Blur Circle:}
	Due to spherical aberration, rays that are not paraxial do not converge to a single point, leading to the 
	formation of a blur circle on the image plane. The radius of this blur circle $r_b$ can be approximated by
	assuming that the deviation in focal points, known as longitudinal spherical aberration, grows quadratically with 
	the ray height $y$ on the mirror surface:
	
    \begin{equation}
    \Delta d \propto y^2
    \end{equation}
	
    Given that the maximum ray height is approximately $D/2$, this suggests:
	
    \begin{equation}
    \Delta d \propto \left(\frac{D}{2}\right)^2 \propto D^2.
    \end{equation}

    The resulting transverse blur (i.e., the blur circle radius) scales with the deviation in the focus 
	position $\Delta d$ and the ratio $D/f_0$, leading to:
	
    \begin{equation}
    r_b \sim \Delta d \cdot \frac{D}{f_0} \propto \frac{D^3}{f_0 R}.
    \end{equation}
	
    This captures the essential reason why the blur circle radius grows with the cube of the aperture diameter.

   \item \textbf{Higher-Order Aberration Terms:}
   The blur circle's radius increases with the fourth power of the mirror's aperture, \( r_b \propto D^4 \), 
   indicating that spherical aberration is a higher-order aberration\footnote{The $D^3$dependency is 
   generally the primary effect seen with spherical aberration, while the $D^4$ dependency reflects 
   higher-order corrections that become significant as the aperture size increases or as more detailed 
   aberration terms are considered in the optical analysis.}. This relationship emphasizes the 
   non-linear nature of the aberration, where small increases in aperture size lead to disproportionately 
   large increases in aberration.

   \item \textbf{Mitigating Spherical Aberration:}
   In optical design, spherical aberration can be minimized by using aspherical mirrors or lenses, which 
   have a curvature that varies with the distance from the optical axis. These designs aim to bring all 
   rays to a common focus, thereby reducing the blur circle and improving image quality.

\end{enumerate}

As a result of these aberrations, the anticipated shape of the Focal Body (FB) resembles a \textbf{mushroom-like} structure, with a denser, 
more concentrated region near the optical axis (the \emph{stem}) and a broader, curved caustic envelope extending outward (the \emph{cap}). 
This morphology is a direct consequence of the gradual deviation of marginal rays, which intersect the optical axis closer to the 
mirror due to spherical aberration.

For a comprehensive treatment of spherical aberration and other optical aberrations, please refer 
to~\cite{born1999}.

\subsection{Flexible Foil Mesh Generation (FFMG)}

The Flexible Foil Mesh Generation (FFMG) method~\cite{Moriya:2025FlexableFoil} constructs high-quality, closed-surface meshes from confined 
$3D$ point clouds using a physically-based simulation of flexible foils. It integrates elasticity-driven deformation, pressure-induced contraction, and adaptive snapping to fixed vertices, ensuring geometric fidelity and computational efficiency.

A deformable triangular mesh $\mathcal{M} = (V,F)$ is initialized around the point cloud using a convex hull approximation. 

To effectively model the dynamic adaptation of the mesh, we employ a force-equilibrium approach, ensuring that the mesh deforms in 
response to applied forces without inertia effects. This method guarantees that the mesh conforms to the point cloud while maintaining 
structural integrity through controlled deformation.

In this formulation, each vertex is treated as a massless node, connected to its neighbors via elastic constraints that regulate 
local stiffness and prevent excessive distortion. Additionally, an external pressure force drives the contraction of the mesh, 
mimicking the behavior of a flexible foil minimizing its surface area. To stabilize the evolution, a numerical damping term is 
introduced, preventing large, unstable steps and ensuring smooth convergence.

The governing equation for vertex motion is given by:

\begin{equation}
    0 = f_{\text{elastic},i} + f_{\text{pressure},i} - c \frac{d u_i}{dt} \label{eq:dynamicElastPress}
\end{equation}

where $u_i$ is the spatial position of the $i$'th mesh-point vertex, $f_{\text{elastic},i}$ represents elastic 
restoration forces, $f_{\text{pressure},i}$ accounts for uniform or spatially varying pressure loads, and $c$ is a damping coefficient 
that stabilizes the simulation.

This ensures a quasi-static evolution rather than a fully dynamic simulation, meaning that the mesh continuously deforms under 
applied forces without oscillatory motion associated with inertia.

Elasticity is modeled through a vertex-based elastic force formulation, where each vertex is influenced by the 
relative displacement of its neighbors:

\begin{equation} 
    f_{\text{elastic},i} = \sum_{j \in N(i)} k_{ij} \left(\frac{\|u_j - u_i\| - L_{ij}}{L_{ij}}\right) \frac{(u_j - u_i)}{\|u_j - u_i\|},
\end{equation}

where $k_{ij}$ is the stiffness coefficient, $L_{ij}$ is the rest length of edge $(i,j)$, and $N(i)$ denotes the neighboring vertices.

Since the model operates in an \textbf{overdamped regime}, the acceleration term is neglected, and the motion of vertices is determined 
purely by the balance of forces as in \eqref{eq:dynamicElastPress} above.

This formulation ensures that vertex displacements occur in response to external forces without inertial effects, meaning the system 
continuously evolves towards a lower-energy configuration without oscillatory motion.

The pressure force on each vertex is computed as:

\begin{equation}
    f_{\text{pressure},i} = \frac{p \cdot n_f \cdot A_f}{3},
\end{equation}

where $p$ is the pressure magnitude exerted on facet $f$, $n_f$ is the normal of the facet, and $A_f$ is its area. 
Fixed vertices $V_{\text{fixed}} \subset V$ constrain the deformation, and a snapping mechanism ensures 
proximity-based vertex alignment.

The mesh evolution follows an explicit Euler integration scheme:

\begin{equation} \label{eq:euler_integration}
    v_i^{t+1} = v_i^t + \Delta t \cdot a_i^t,
\end{equation}
\begin{equation}
    x_i^{t+1} = x_i^t + \Delta t \cdot v_i^{t+1},
\end{equation}

where $\Delta t$ is the time step and $v_i^t$ and $a_i^t$ are the velocity and respective acceleration computed from net forces.

To prevent numerical instability, the method enforces a Courant-Friedrichs-Lewy (CFL) condition \cite{CLF1928}:

\begin{equation}
    \Delta t < \frac{2}{\omega_{\max}},
\end{equation}

where $\omega_{\max}$ is the highest eigenfrequency of the system.

By iteratively refining the mesh structure through adaptive smoothing and convergence testing, FFMG generates a physically realistic, topologically consistent mesh that conforms to complex point cloud geometries.

\paragraph{Minimal Surface Approximation}
The physically-based simulation inherently favors minimal-area solutions within the constraints imposed by the 
fixed vertices. The pressure-driven contraction mimics the behavior of a flexible membrane seeking an 
energy-minimizing configuration, akin to soap films forming minimal surfaces under boundary constraints. 
The elastic forces work to maintain structural integrity while minimizing local surface tension, resulting in a 
final mesh that approximates a constrained minimal surface. While the method does not explicitly solve a
minimal surface equation, its dynamic evolution naturally leads to a stable, near-minimal surface 
configuration, subject to the imposed geometric constraints.

\section{Methodology}

\subsection{Problem Definition}

In this study, we seek to accurately determine the $3D$ Focal Body generated by a spherical mirror 
focusing light rays from an infinitely distant source along its optical axis. 
We use a ray-tracing simulation to follow a set of rays reflected by the mirror and count their
mutual intersections as a function of occurrence in $3D$ space.
The Focal Body is the region where the reflected rays converge, forming a complex, highly concave structure. 
Traditional methods struggle to capture the full extent of this concave geometry, which is crucial 
for precise optical modeling.

Formally, let \( P = \{p_1, p_2, \dots, p_n\} \) represent the set of points in $3D$ space that describe 
the confined volume of the Focal Body. These points are derived from the intersection of reflected light rays 
with a hypothetical observation plane placed near the focal region. The goal is to construct a 
surface \( S^{cc} \) that:

\begin{enumerate}
    \item Encloses all points in \( P \), ensuring no points are excluded.
    \item Tend to minimizes the surface area while maintaining a concave structure that accurately represents the 
	complex geometry of the Focal Body.
    \item Remains closed and non-intersecting to preserve the integrity of the surface, ensuring that the 
	Focal Body is fully encapsulated without self-intersections.
\end{enumerate}

The FFMG method is employed to construct a surface \( S^{cc} \) that meets the 
above criteria by iteratively refining an initial convex hull. The process involves replacing facets on the 
convex hull with new facets that better approximate the concave geometry of the Focal Body.

Mathematically, the problem can be formulated as finding a surface \( S^{cc} \) that approximately minimizes 
surface area while satisfying the necessary constraints:

\begin{equation}
S^{cc} \approx \argmin_{S} \left\{ \text{Area}(S) \mid P \subseteq \text{Int}(S), \, S \text{ is non-intersecting and closed} \right\}
\end{equation}

where \( \text{Area}(S) \) represents the surface area of \( S \), and \( \text{Int}(S) \) denotes the 
interior region enclosed by \( S \). 

However, the FFMG method does not explicitly solve a minimal surface problem; rather, it approximates a 
constrained minimal surface through a physically-based simulation. The method relies on pressure-driven 
contraction and elastic forces to evolve the surface dynamically, tending toward a locally minimal 
configuration within the constraints imposed by fixed vertices. While this process favors a minimal-area solution, 
the final surface is shaped by both physical simulation and numerical constraints, leading to an approximate 
rather than an exact minimization of surface area.

In summary, the problem involves using FFMG to generate a $3D$ surface that closely captures the concave nature 
of the light convergence in a spherical mirror system. This provides a detailed yet computationally efficient 
approximation essential for optical design and analysis.

\subsection{Ray Tracing Computational Implementation}

The Ray Tracing simulation is developed to model and analyze the reflection behavior of light rays 
interacting with complex, arbitrary-shaped reflecting mirror geometries. 

\vspace{0.5cm}

\begin{figure}[H]
\begin{center}
\includegraphics[width=7cm]{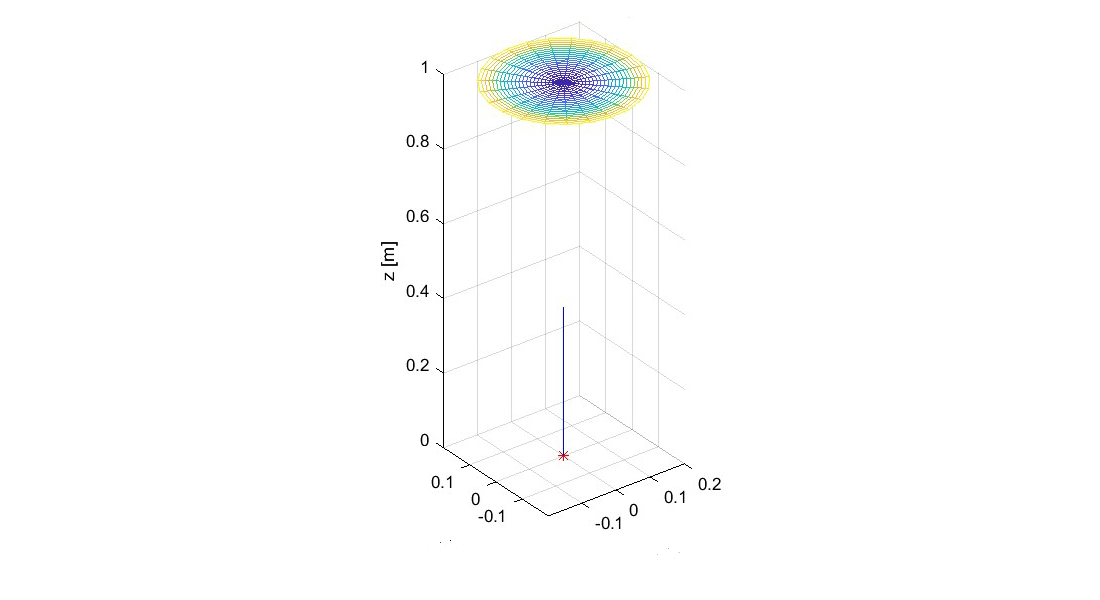}
\includegraphics[width=7cm]{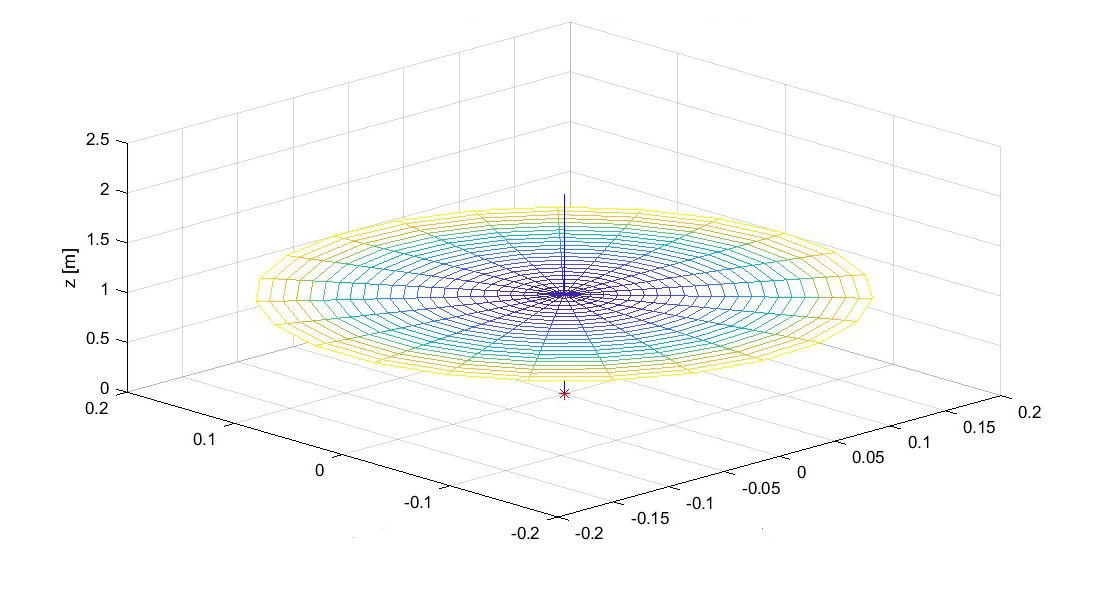}
\caption{A typical simulation setup for spherical mirror with an equal-axis view (left) and non-equal 
view (right).}
\label{A typical simulation setup for spherical mirror}
\end{center}
\end{figure}

The key features of the simulation include:

\begin{itemize}

\item \textbf{Mirror Geometry}: The simulation primarily focuses on cylindrical mirrors, including 
configurations with segmented, arbitrarily-shaped mirrors exhibiting cylindrical symmetry. 
This approach enables the exploration of a wide range of setups, from single mirrors to complex 
arrays of multiple cylindrical mirrors (see in fig. \ref{A typical simulation setup for spherical mirror}).

\item \textbf{Efficiency through Simplification}: To maintain computational efficiency, the simulation 
employs cylindrical mirrors to approximate non-cylindrical symmetric mirrors. 
Logical functions are used to determine ray intersections, effectively reducing the simulation's 
complexity.

\item \textbf{Optimization and Consistency}: Optimization is applied not at the level of individual 
mirrors but to the overall system configuration, ensuring consistency with established results for 
arbitrarily-shaped mirrors and facilitating smooth transitions between different setups.

\item \textbf{Ray Source and Reflection Dynamics}: The simulation handles multiple light sources, 
including a primary on-axis source and additional off-axis sources. The ray density is adjusted to 
account for losses due to mirror segmentation, ensuring accurate reflection modeling.

\item \textbf{Comprehensive Analysis of Results}: The simulation offers an in-depth analysis of 
various mirror configurations, including single spherical mirrors, segmented parabolic mirrors, 
arbitrarily-shaped mirrors, and complex arrays composed of diverse sets of smaller mirrors. 
This analysis provides valuable insights into how different segmentations and orientations 
impact the distribution of reflected light.

\item \textbf{Optimization Techniques}: The simulation includes robust optimization algorithms to 
adjust the position and orientation of the mirrors, considering various degrees of freedom such as 
pitch, azimuth, and spatial coordinates, to minimize image distortion and optimize focus.

\item \textbf{Source and Ray Parameters}: The simulation allows for the configuration of multiple 
light sources, each with a customizable number of rays. These rays are traced from their sources 
to the mirrors and ultimately to a planar surface (analogous to a CCD), where the resulting image 
quality is thoroughly analyzed.

\item \textbf{Numerical and Analytical Methods}: The simulation employs both numerical and analytical 
methods for computing ray reflections, ensuring the accuracy and reliability of the results by validating 
the methods against each other.

\item \textbf{Output and Analysis}: Detailed outputs, including images and data files, are generated, 
capturing key metrics such as spot sizes and image quality for different configurations. 
These outputs are crucial for evaluating the effectiveness of various mirror setups and optimization 
strategies.

\item \textbf{planar surface and Focus Optimization}: A specialized module is included for optimizing the  
planar surface position to minimize spot size, involving calculations of standard deviation (STD) or root 
mean square (RMS) focus, and adjusting planar surface height to achieve optimal focus.

\item \textbf{Reflected Rays}: The simulation produces an output matrix that precisely tracks each 
incoming ray, detailing its interaction with the mirror. This includes the coordinates of the 
impact point on the reflecting surface and the orientation of the reflected ray in space, all as 
functions of the specific mirror segment the ray encounters.

\item \textbf{Ray Intersections}: In the post-simulation analysis, the reflected rays are evaluated 
against a predefined $3D$ grid. The simulation records and creates histograms of the number of rays 
passing through each cell in this grid, providing a detailed spatial distribution of ray intersections.

\item \textbf{Focal Body Definition}: Utilizing histogrammed layers based on the density of rays (from 
the highest value down to a predefined minimum, typically greater than 1\% of the total incoming rays), 
a shell-based Focal Body is constructed. Each layer defines a bounded volume encompassing all cells 
through which rays of that histogram value or higher pass. At this stage, the Focal Body is represented 
as a point cloud.

\end{itemize}

\subsection{Application of FFMG to the Focal Body}
This section describes the step-by-step application of the Flexible Foil Mesh Generation (FFMG) 
method to simulate the $3D$ Focal Body formed by a spherical mirror. The Focal Body is a concave 
structure resulting from the convergence of light rays reflected by the mirror. The FFMG method 
generates a mesh that approximates this shape by iteratively refining an initial convex approximation 
to better conform to the concave geometry of the Focal Body.

\subsubsection{Step 1: Initialization with Convex Hull}

To ensure consistency in scale, the point cloud is normalized such that each dimension is linearly mapped to the range \([-1,1]\). 
This transformation is defined as:

\begin{equation}
    \tilde{p}_i = \frac{2 (p_i - p_{\min})}{p_{\max} - p_{\min}} - 1,
\end{equation}

where \( p_{\min} \) and \( p_{\max} \) denote the minimum and maximum values along each dimension, respectively. 
As a result, all points satisfy:

\begin{equation}
    \tilde{p}_i \in [-1,1]^d,
\end{equation}

ensuring a standardized representation for mesh reconstruction.

Following the normalization stage, the convex hull \( S^{ch} \) of the point cloud \( P \)is computed. 
This convex hull serves as the initial approximation of the surface, representing the starting boundary of the Focal Body:

\begin{equation}
S^{ch} = \text{Conv}(P) = \bigcap \{ H \mid H \text{ is convex, } P \subseteq H \}.
\end{equation}

While this surface encloses all points in \( P \), it does not capture the concave regions of the 
Focal Body and requires further refinement.

\subsubsection{Step 2: Convex Hull Refinement}
After computing the convex hull \( S^{ch} \), the next step involves refining the mesh by introducing 
denser triangular facets, particularly in regions near the convex hull vertices. The goal of this refinement 
is to increase the resolution of the initial approximation before applying further deformation.

The refinement process subdivides the triangular facets of \( S^{ch} \) by introducing new vertices along 
edges and within selected facets. This is achieved through a recursive subdivision method that maintains 
the integrity of the convex hull while ensuring a more uniform and adaptable triangulation.

A key aspect of this refinement process is the adaptive density increase near convex hull vertices. 
Let \( v_i \) be a vertex in \( S^{ch} \), and let \( N(i) \) be the set of neighboring vertices. 
To determine whether refinement is needed, a local density criterion is applied:

\begin{equation}
L_i = \frac{1}{|N(i)|} \sum_{j \in N(i)} \| v_j - v_i \|,
\end{equation}

where \( L_i \) represents the average edge length surrounding \( v_i \). A vertex \( v_i \) undergoes 
further subdivision if:

\begin{equation}
L_i > \gamma L_{\text{min}},
\end{equation}

where \( \gamma > 1 \) is a refinement factor, and \( L_{\text{min}} \) is the target edge length in 
high-resolution regions.

Each refined facet is subdivided by inserting midpoints along its edges, generating new triangular facets. 
For a given triangle \( (v_i, v_j, v_k) \), the midpoints \( m_{ij}, m_{jk}, m_{ki} \) are computed as:

\begin{equation}
m_{ij} = \frac{v_i + v_j}{2}, \quad m_{jk} = \frac{v_j + v_k}{2}, \quad m_{ki} = \frac{v_k + v_i}{2}.
\end{equation}

These midpoints replace the original triangle with four new sub-triangles, maintaining a non-intersecting 
and closed structure. The subdivision process continues iteratively until the average edge length reaches the 
desired refinement level.

The refinement stage ensures that the initial convex hull provides a sufficiently detailed starting point for 
subsequent deformation, allowing the mesh to better conform to concave regions in later stages.

\subsubsection{Step 3: Iterative Surface Refinement}
Once the initial set of facets has been constructed, the mesh undergoes an iterative refinement process 
driven by physically-based deformation. The mesh \( S^{cc}_t \) at iteration \( t \) evolves 
according to the force balance equation \eqref{eq:dynamicElastPress}.

The vertices positions \( u_i \) are updated using explicit time integration equation \eqref{eq:euler_integration}
and the refinement process continues until a convergence criterion is met:

\begin{equation}
\max_{u_i \in V} \| u_i^{t+1} - u_i^t \| < \epsilon,
\end{equation}

where \( \epsilon \) is a predefined threshold ensuring that the surface has reached a stable 
configuration.

\subsubsection{Step 4: Modifications and Optimizations}
Several modifications and optimizations are applied to improve the quality and stability 
of the generated surface:

- \textit{Geometric smoothing:} A Laplacian smoothing operator is applied to prevent excessive 
local deformations:

\begin{equation}
u_i' = u_i + \lambda \sum_{j \in N(i)} (u_j - u_i),
\end{equation}

where \( \lambda \) is a smoothing parameter, and \( N(i) \) represents the neighboring vertices.

- \textit{Adaptive snapping:} Certain mesh-points vertices are constrained to snap onto predefined fixed positions 
(the point cloud vertices) to align the mesh with the reference geometry:

\begin{equation}
u_i' = \argmin_{u_j \in V_{\text{fixed}}} \| u_i - u_j \|.
\end{equation}

- \textit{Topology preservation:} The refinement process maintains a non-self-intersecting constraint 
to prevent topological inconsistencies.

\subsubsection{Step 5: Final Surface Validation}
After the iterative refinement process, the final surface \( S^{cc} \) is validated to ensure that 
it meets the expected geometric and topological criteria. This validation consists of:

- \textit{Enclosure test:} Verifying that all points in \( P \) satisfy

\begin{equation}
P \subseteq \text{Int}(S^{cc}).
\end{equation}

- \textit{Self-intersection check:} Ensuring that the surface remains free of self-intersections.

- \textit{Curvature analysis:} Evaluating whether the surface curvature is consistent with the 
expected concave structure of the Focal Body.

The application of FFMG to the $3D$ Focal Body involves an iterative process of surface refinement, 
which transforms an initial convex hull into a concave representation of the Focal Body. The 
refinement steps, combined with geometric constraints and stability enhancements, ensure that the 
generated mesh closely approximates the underlying concave geometry while maintaining computational 
efficiency.

\subsection{Computational Complexity Analysis}

The Flexible Foil Mesh Generation (FFMG) method involves a series of computationally intensive steps, each contributing to the overall complexity of the simulation. These steps include mesh initialization, force computation, iterative deformation updates, and convergence testing. Below, we present a high-level analysis of the computational complexity associated with the major components of the simulation.

\subsubsection{Mesh Initialization and Subdivision}
The initial mesh is generated by computing the Convex Hull of the input point cloud, which requires \(\mathcal{O}(N \log N)\) operations using the QuickHull algorithm. This provides an initial enclosing surface for the Focal Body. To improve resolution, a subdivision process is applied to the mesh, where each triangular face is iteratively divided, increasing the number of elements by a factor of \(4^s\), where \(s\) is the subdivision level. This results in an exponential increase in the number of faces, leading to an approximate computational complexity of \(\mathcal{O}(M)\), where \(M\) denotes the number of faces after subdivision.

\subsubsection{Force Computations and Deformation Updates}
Each simulation iteration involves computing two primary forces:
\begin{enumerate}
    \item \textbf{Pressure Forces}: These are computed per face, leading to a complexity of \(\mathcal{O}(M)\).
    \item \textbf{Elastic Forces}: Modeled using a force-equilibrium system, each vertex interacts with a fixed number of neighboring vertices, 
	leading to a complexity of \(\mathcal{O}(N)\).
\end{enumerate}
Additionally, vertex displacement updates occur using explicit time integration, which scales as \(\mathcal{O}(N)\). Since the number of faces and vertices are related by \(M \approx 2N\) for typical triangular meshes, the force computation and deformation update complexity is approximated as \(\mathcal{O}(N)\) per iteration.

\subsubsection{Nearest-Neighbor Searches and Snapping}
The snapping mechanism ensures that vertices near fixed points are properly aligned during deformation. This requires repeated nearest-neighbor searches, implemented efficiently using a KD-tree with \(\mathcal{O}(N \log N)\) complexity for tree construction and \(\mathcal{O}(\log N)\) per query. Given that nearest-neighbor searches are performed for all \(N\) vertices, this results in an overall complexity of \(\mathcal{O}(N \log N)\) per iteration.

\subsubsection{Convergence and Iteration Scaling}
The simulation progresses over multiple iterations until convergence criteria are met. In each iteration, the maximum vertex displacement is 
checked against a predefined tolerance, requiring \(\mathcal{O}(N)\) operations. 

In this study, a fixed time-stepping approach is used to update the mesh deformation process. 
However, for complex concave structures, adaptive time-stepping could offer an alternative strategy to 
dynamically adjust step size based on deformation rates, improving stability and efficiency.

A key tradeoff exists between fixed and adaptive time-stepping: while fixed time steps ensure predictable iteration 
counts and easier convergence analysis, they may require excessively small steps to maintain numerical stability, 
leading to increased computational costs. Conversely, adaptive time-stepping dynamically adjusts $\Delta t$
based on local deformation magnitudes, potentially accelerating convergence for smoothly evolving regions while preventing 
overshooting in high-curvature areas.

The total number of iterations is bounded, typically in the range \(50 \leq k \leq 500\). Thus, the total computational cost of the 
simulation is approximately:

\[
\mathcal{O}(k (N \log N + N)) = \mathcal{O}(k N \log N).
\]

\subsubsection{Overall Complexity Estimate}
Summarizing the dominant computational steps:

\begin{table}[h]
    \centering
    \begin{tabular}{|l|c|}
        \hline
        \textbf{Computation Step} & \textbf{Complexity} \\
        \hline
        Convex Hull Initialization & \(\mathcal{O}(N \log N)\) \\
        Mesh Subdivision & \(\mathcal{O}(M)\) \\
        Force Computation & \(\mathcal{O}(N)\) \\
        Nearest-Neighbor Snapping & \(\mathcal{O}(N \log N)\) \\
        Iterative Updates & \(\mathcal{O}(k N \log N)\) \\
        \hline
    \end{tabular}
    \caption{Computational complexity of major steps in the simulation.}
\end{table}

Since the number of faces after subdivision grows with \(M = \mathcal{O}(4^s N)\), choosing a high subdivision level significantly increases the computational cost. However, the iterative updates remain the dominant term, making the simulation complexity approximately \(\mathcal{O}(k N \log N)\) in practical scenarios.

\subsubsection{Optimization Considerations}
To improve computational efficiency, several optimizations can be applied:
\begin{itemize}
    \item \textbf{Reducing subdivision levels} limits exponential mesh growth.
    \item \textbf{Parallelizing force computations} using vectorized operations can reduce per-iteration costs.
    \item \textbf{Sparse data structures} can be used to minimize unnecessary computations on inactive regions of the mesh.
    \item \textbf{Adaptive time-stepping} can reduce the number of required iterations by dynamically adjusting the deformation step size.
\end{itemize}

The computational complexity of the FFMG simulation is primarily governed by the number of vertices, faces, and the number of iterations required for convergence. The nearest-neighbor search and force computations dominate each iteration, leading to an overall complexity of \(\mathcal{O}(k N \log N)\). The framework remains computationally feasible for large-scale point clouds, but careful selection of parameters such as subdivision level and convergence tolerance is necessary to ensure efficiency.

\section{Results}

\subsection{Simulation Setup}

The simulation was designed to model the deformation of a flexible foil under the influence of pressure forces and elastic constraints. 
The following parameters were employed to ensure numerical stability, realistic physical behavior, and computational efficiency.

\subsubsection{Mesh Subdivision and Initialization}
\begin{itemize}
    \item \textbf{subdivision\_level = 3}: This parameter defines the initial refinement level of the mesh, with higher values leading to a finer discretization of the foil. A subdivision level of 3 ensures a balance between computational efficiency and geometric accuracy. Typically, values range from \textbf{1} (coarse) to \textbf{5} (very fine), where higher levels improve resolution but increase computational cost.
\end{itemize}

\subsubsection{Pressure and External Forces}
\begin{itemize}
    \item \textbf{PR\_in = 0.5, PR\_out = 1.0}: These parameters represent the internal and external pressure applied to the foil, respectively. The net pressure difference, defined as \( P = PR_{\text{in}} - PR_{\text{out}} \), determines whether the structure expands or contracts. In this setup, \( P = -0.5 \) induces a compressive effect on the foil.
    \item \textbf{pressure\_scaling\_factor = 10.0}: This factor amplifies the magnitude of pressure forces to ensure they have a significant 
	impact on deformation. Typical values range from \textbf{1 to 50}, depending on the scale of the model and the physical properties of the material.
    \item \textbf{Pressure force mode = globalCoM}: Specifies the methodology used to compute pressure force directions. The \textbf{global center of mass (globalCoM) method} ensures forces act consistently relative to the centroid of the fixed vertices. Alternative modes, such as \textit{Norm} (face-normal based) or \textit{localCoM} (local neighborhood-based), may be used for different physical behaviors.
\end{itemize}

\subsubsection{Incremental Pressure Application}
\begin{itemize}
    \item \textbf{pressure\_increment = 0.03}: Defines the incremental increase in applied pressure per iteration. This gradual application prevents numerical instability. Values are typically \textbf{0.01–0.05}, depending on the required resolution of deformation tracking.
\end{itemize}

\subsubsection{Snapping and Mesh Constraints}
\begin{itemize}
    \item \textbf{snapping\_tolerance = 0.02}: Controls the maximum allowable distance for a mesh vertex to be snapped to a fixed vertex. Ensuring snapping within this tolerance prevents artificial stretching or discontinuities. Values range from \textbf{0.001} (strict) to \textbf{0.05} (loose).
    \item \textbf{max\_NNsnapping\_iterations = 5}: Specifies the number of iterations allowed for nearest-neighbor snapping relaxation, preventing excessive distortions due to sudden vertex movements.
\end{itemize}

\subsubsection{Deformation Criteria and Convergence}
\begin{itemize}
    \item \textbf{deformation\_tolerance = 1e-5}: Defines the \textbf{convergence criterion} for the simulation. The deformation process halts when the maximum displacement of all vertices falls below this threshold. Values typically range from \(10^{-6}\) (very strict) to \(10^{-3}\) (loose).
    \item \textbf{mTol = 0.016}: This parameter constrains the \textbf{maximum allowable displacement} per iteration to prevent excessive motion. It is set to 80\% of \textit{snapping\_tolerance}, ensuring controlled deformation.
\end{itemize}

\subsubsection{Temporal Resolution and Iterative Processing}
\begin{itemize}
    \item \textbf{deformation\_max\_iterations = 200}: The \textbf{maximum number of deformation iterations} before termination. Higher values allow more gradual convergence but increase computational time. This is typically set between \textbf{50 and 500}, depending on simulation complexity.
    \item \textbf{dt = 0.03}: The \textbf{time step} used for numerical integration. It determines how much the simulation progresses per iteration. A smaller \( dt \) results in finer temporal resolution but requires more iterations. Typical values are \textbf{0.01–0.05}.
\end{itemize}

\subsubsection{Smoothing and Damping Mechanisms}
\begin{itemize}
    \item \textbf{smoothingIterations = 1}: Defines the number of \textbf{Laplacian smoothing} passes applied to the mesh to prevent noise artifacts. A value of \textbf{0} disables smoothing. When used, it typically ranges from \textbf{1 to 5}.
    \item \textbf{smoothingTol = 0.02}: Controls the extent of \textbf{smoothing influence}, ensuring excessive mesh modification does not occur. It is linked to \textit{snapping\_tolerance}.
    \item \textbf{damping\_factor = 1.0}: Regulates \textbf{force dissipation}, reducing oscillations in vertex movement. A value of \textbf{1.0} means no damping, while values below \textbf{0.5–0.8} introduce significant damping to stabilize simulations.
\end{itemize}

\subsubsection{Elasticity and Structural Constraints}
\begin{itemize}
    \item \textbf{apply\_snapping = True}: Enables the snapping mechanism, ensuring vertices within \textit{snapping\_tolerance} are relocated to the nearest fixed vertex, preventing unintended distortions.
    \item \textbf{stiffness = 0.01}: Controls the \textbf{elastic resistance} of the foil. A value of \textbf{0} removes elasticity, while higher values (e.g., \textbf{0.05–0.2}) make the foil more resistant to deformation.
    \item \textbf{strain\_factor = 10}: Modulates the \textbf{stiffness response to stretching}. Increasing this value amplifies stiffness at larger strains. Values typically range from \textbf{5 to 15}, depending on material behavior.
    \item \textbf{max\_strain = 0.7}: Specifies the \textbf{maximum allowable strain} before additional stiffness adjustments occur. Values in the range of \textbf{0.5–1.0} provide stability without over-restricting deformation.
    \item \textbf{distance\_factor\_strength = 1}: Controls the \textbf{rate at which fixed vertex influence decays with distance}. Lower values expand the region affected by fixed vertices, while higher values localize their effect.
\end{itemize}

\subsection{Mesh Deformation Results}

\begin{figure}[H]
\begin{center}
\includegraphics[width=11cm]{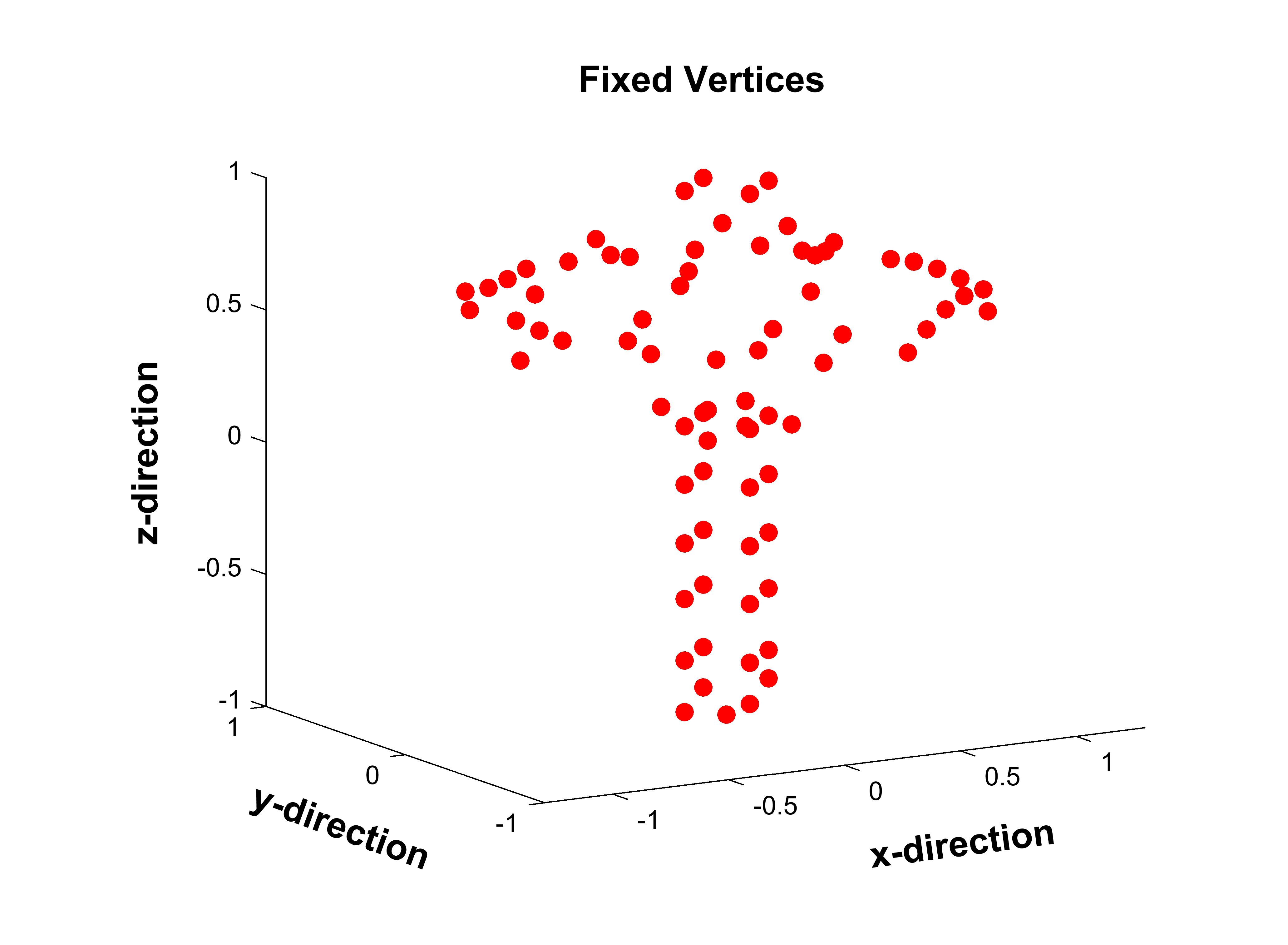}
\caption{Initial point cloud representing the focal body.}
\label{fig:The_Point_Cloud}
\end{center}
\end{figure}

Figure \ref{fig:The_Point_Cloud} illustrates the initial, normalized point cloud utilized in this study, representing the spatial distribution of the focal body. 
The concave structure is evident.
The objective of the mesh deformation process is to generate a closed-surface mesh that accurately conforms to the shape defined by the point cloud.

\begin{figure}[H]
\begin{center}
\includegraphics[width=10cm]{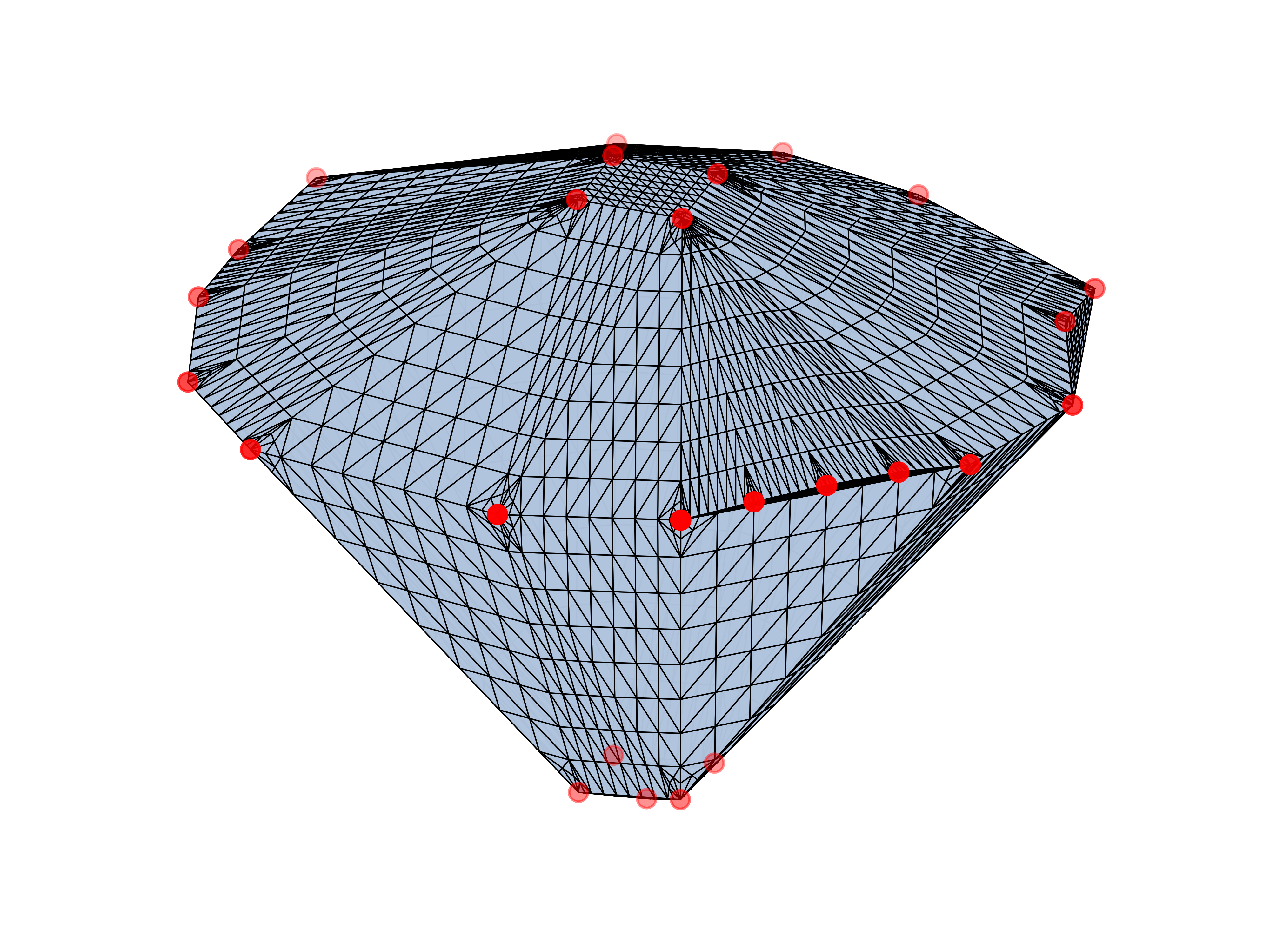}
\caption{Mesh initialization including Convex Hull vertices and refined mesh in their proximity.}
\label{fig:MeshCH}
\end{center}
\end{figure}

The simulation was initiated using the convex hull as the initial approximation of the mesh \ref{fig:MeshCH}. 

\begin{figure}[H]
\begin{center}
\includegraphics[width=7cm]{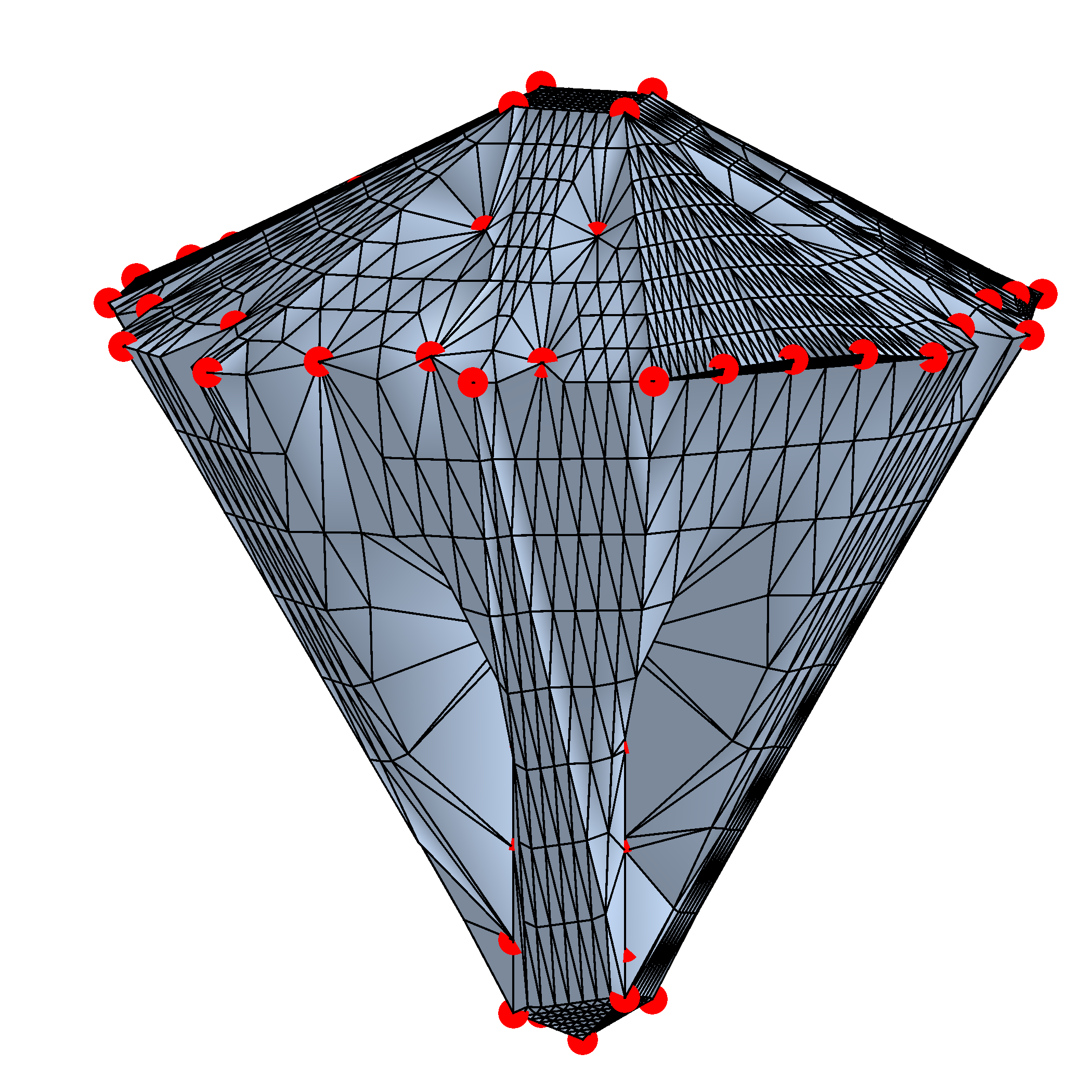}
\includegraphics[width=7cm]{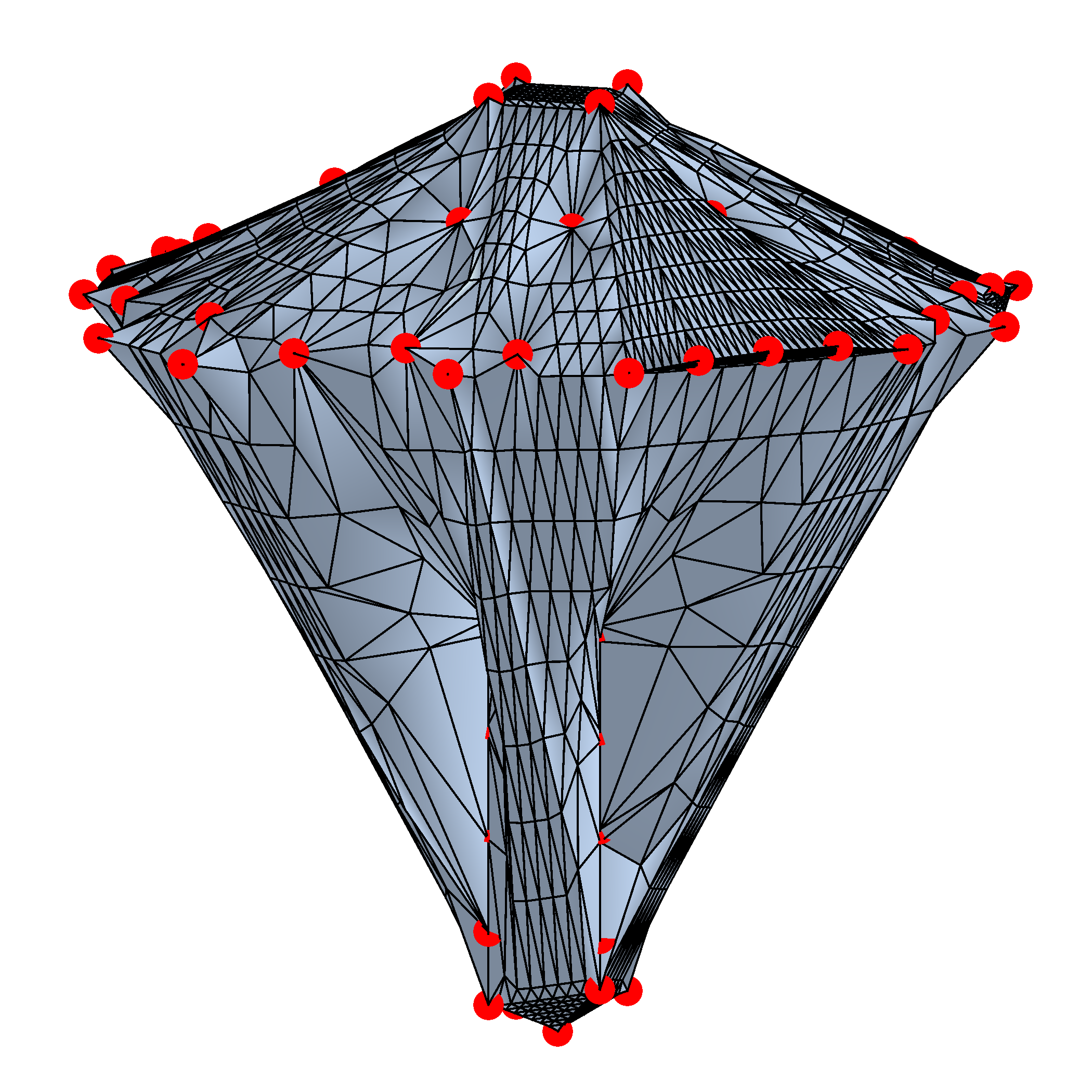} \\
\includegraphics[width=7cm]{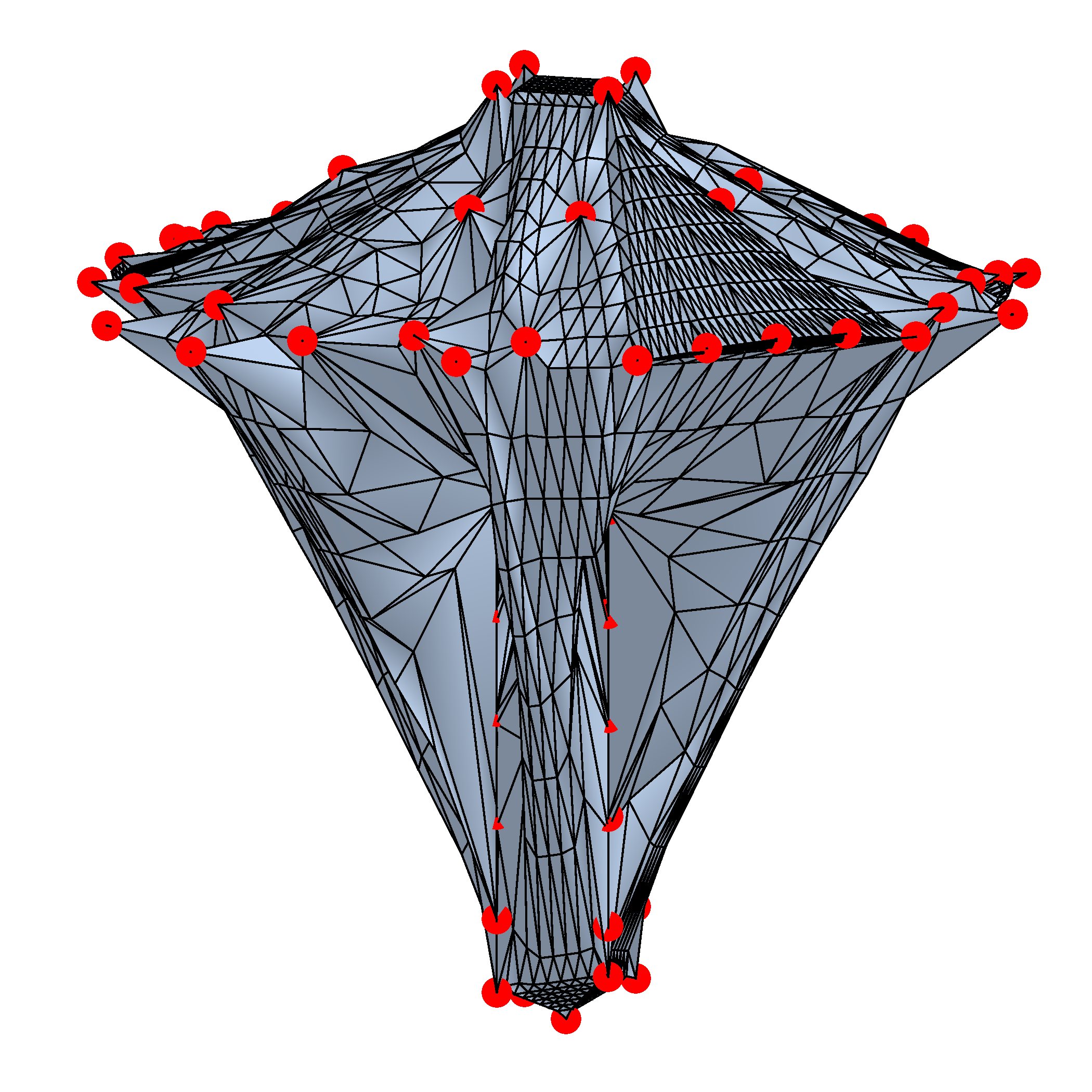}
\includegraphics[width=7cm]{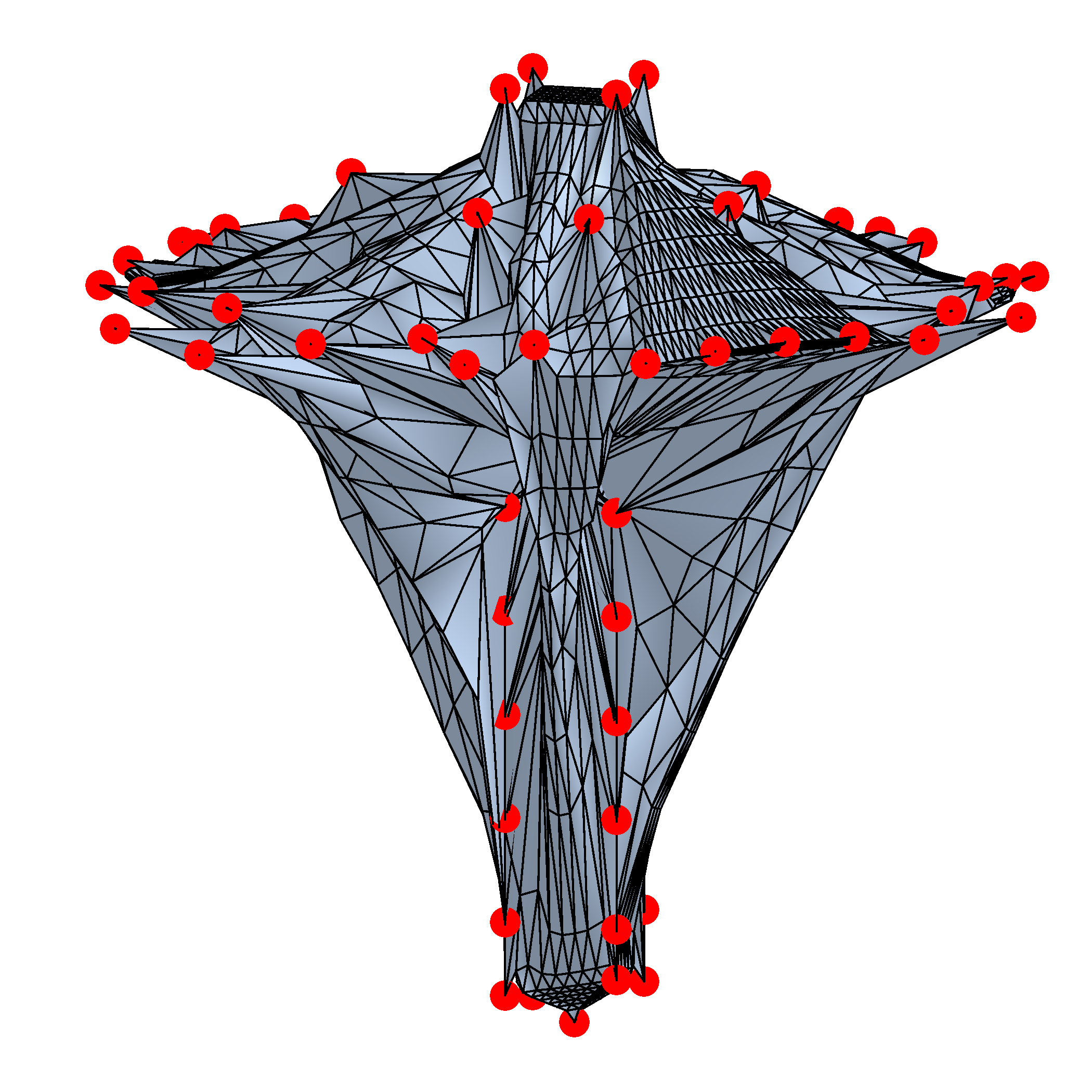}
\caption{Mesh deformation progression at 50 (upper-left), 100 (upper-right), 150 (lower-left), and 199 (lower-right) iterations.}
\label{fig:SimIter}
\end{center}
\end{figure}

Throughout the iterative process, the mesh undergoes deformation to progressively adhere to the underlying point 
cloud structure. Figure \ref{fig:SimIter} illustrates the evolution of the mesh at selected  iterations 50, 100, 150, and 199.

For additional visualization, refer to the following video demonstration, \href{run:Frozen_clip_24Feb2025_180812.mp4}{\textbf{click here}} to 
watch the simulation video.

The results demonstrate that the evolving mesh conforms to the shape of the point cloud while maintaining closure, ensuring a smooth and 
topologically consistent representation of the focal body.

\subsection{Test Point Analysis}

To further analyze the mesh deformation process, specific test points were selected across the mesh grid. The forces acting on these points were tracked over the course of the iterations to gain insights into the dynamics of mesh evolution.

\begin{figure}[H]
\begin{center}
\includegraphics[width=4.5cm]{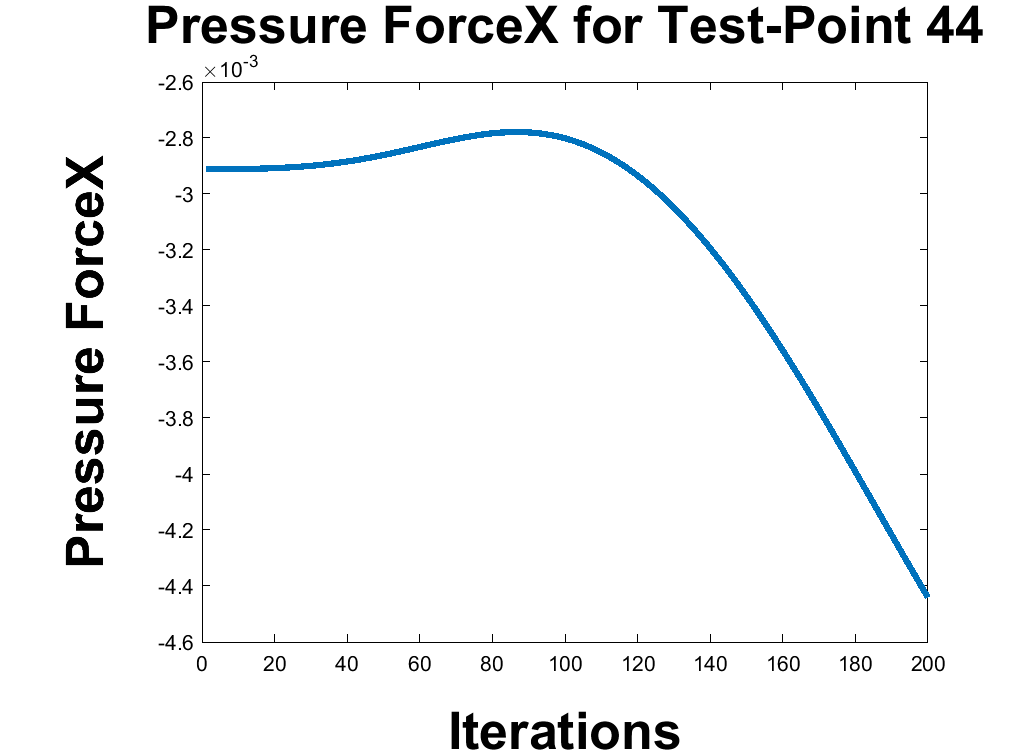}
\includegraphics[width=4.5cm]{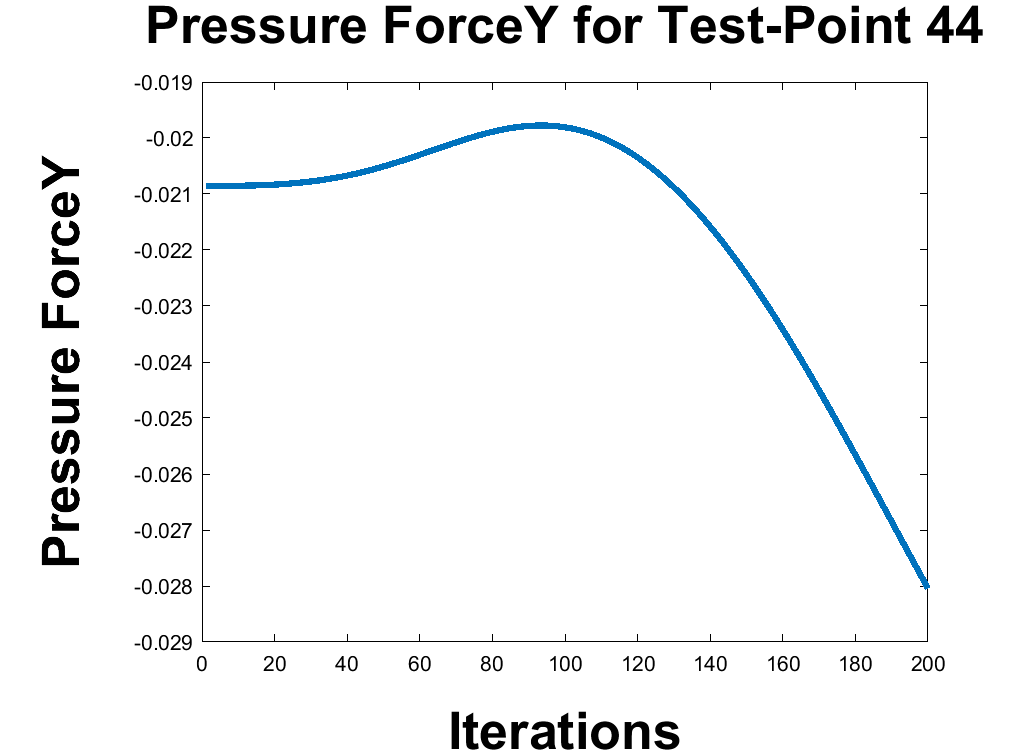}
\includegraphics[width=4.5cm]{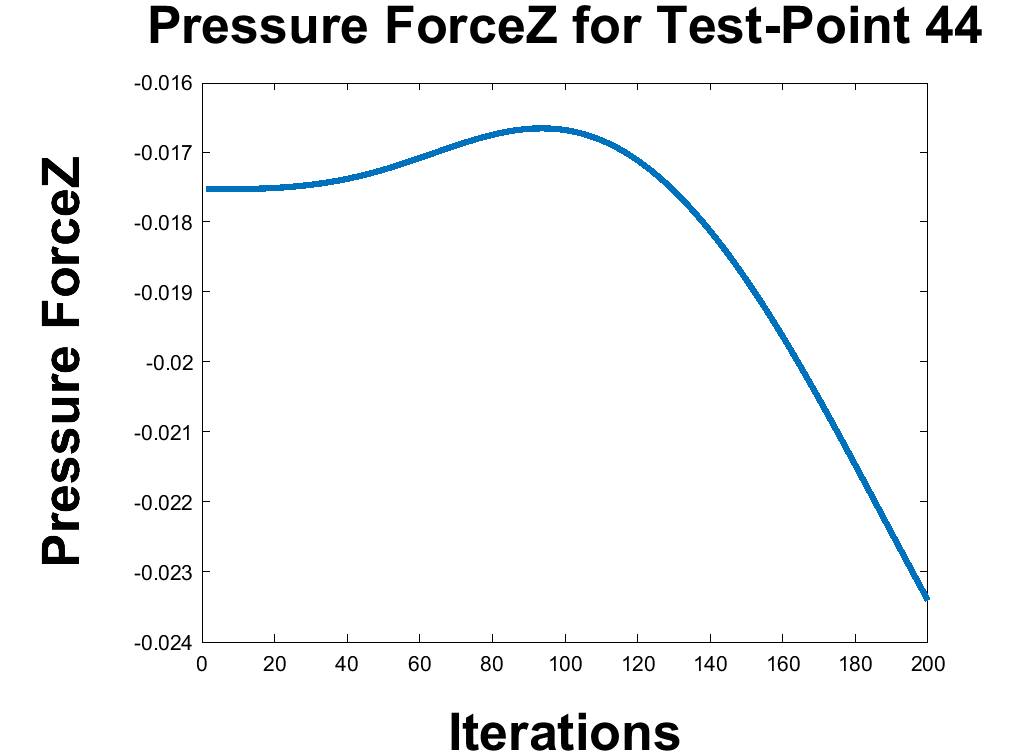}
\caption{Pressure force components acting on a representative test point during the deformation process.}
\label{fig:PressureForceXYZ}
\end{center}
\end{figure}

Figure \ref{fig:PressureForceXYZ} presents the spatial components of the pressure force acting on a representative test point as a function of the simulation iterations.

\begin{figure}[H]
\begin{center}
\includegraphics[width=4.5cm]{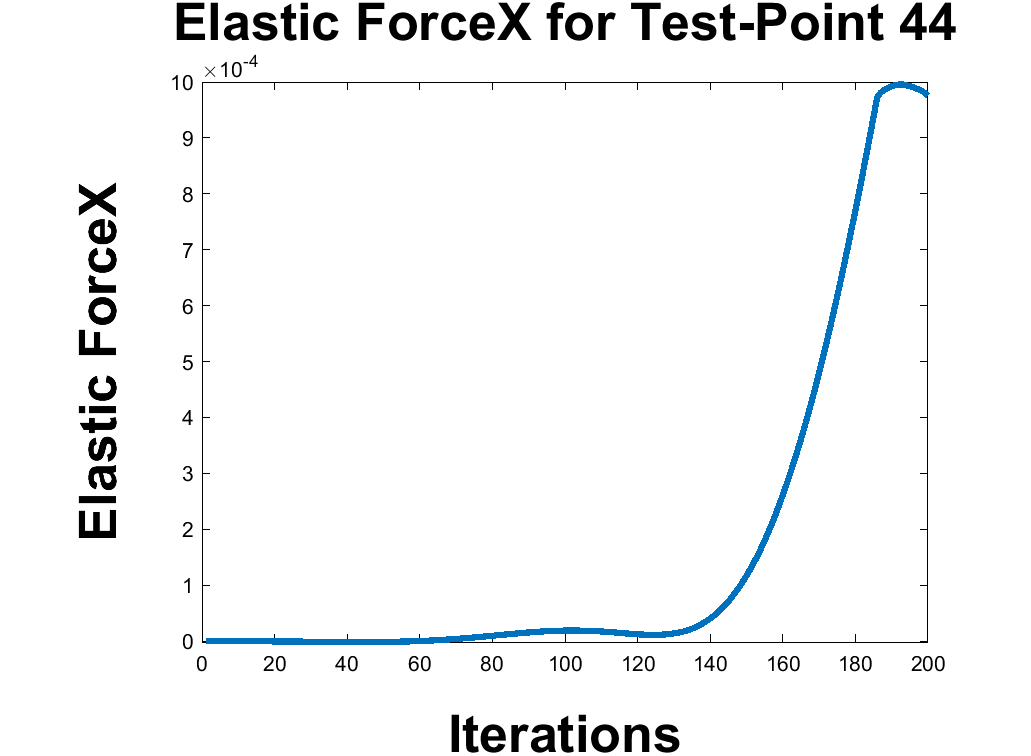}
\includegraphics[width=4.5cm]{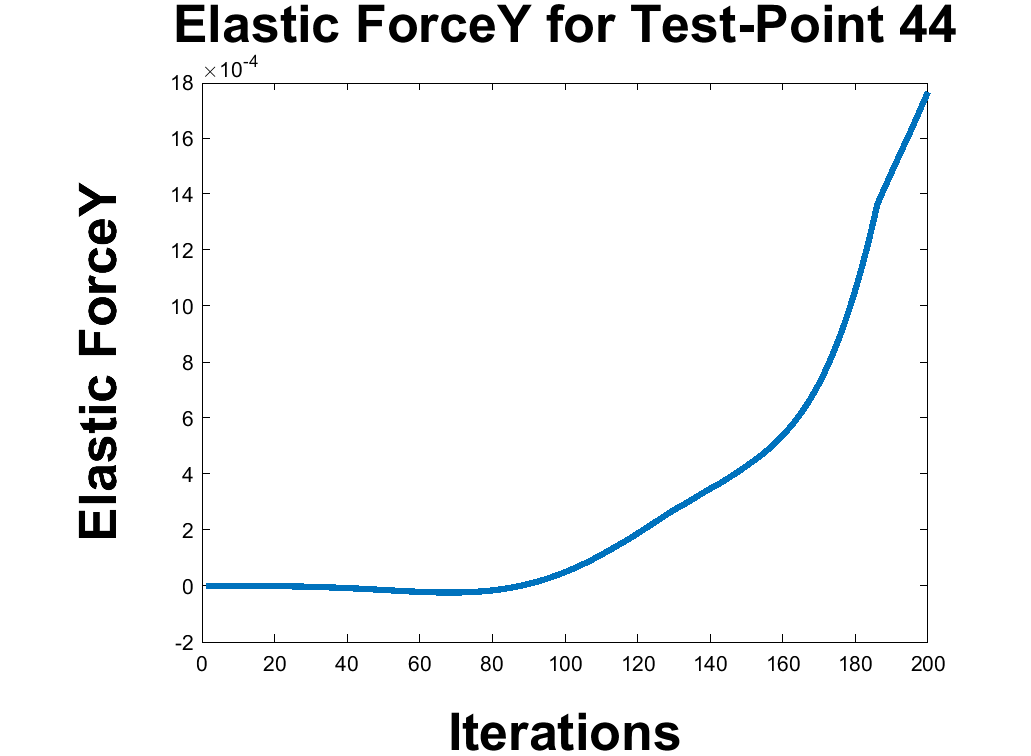}
\includegraphics[width=4.5cm]{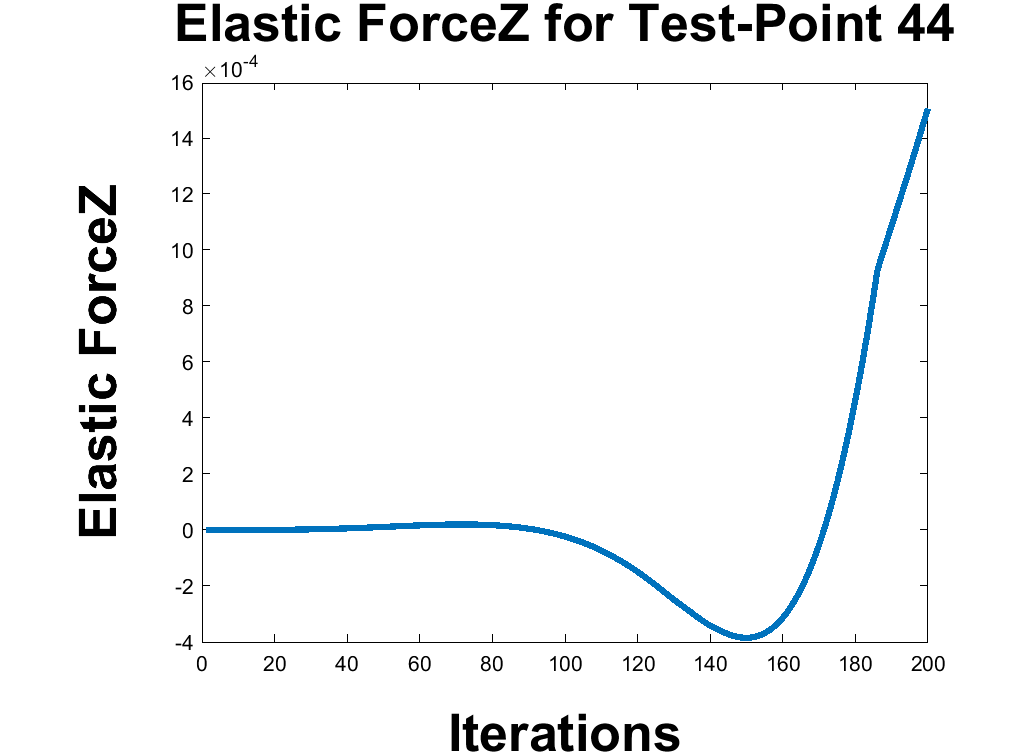}
\caption{Elastic force components acting on the same test point during the deformation process.}
\label{fig:ElasticForceXYZ}
\end{center}
\end{figure}

Similarly, the corresponding elastic force components exerted on the same test point are illustrated in Figure \ref{fig:ElasticForceXYZ}.

\begin{figure}[H]
\begin{center}
\includegraphics[width=7cm]{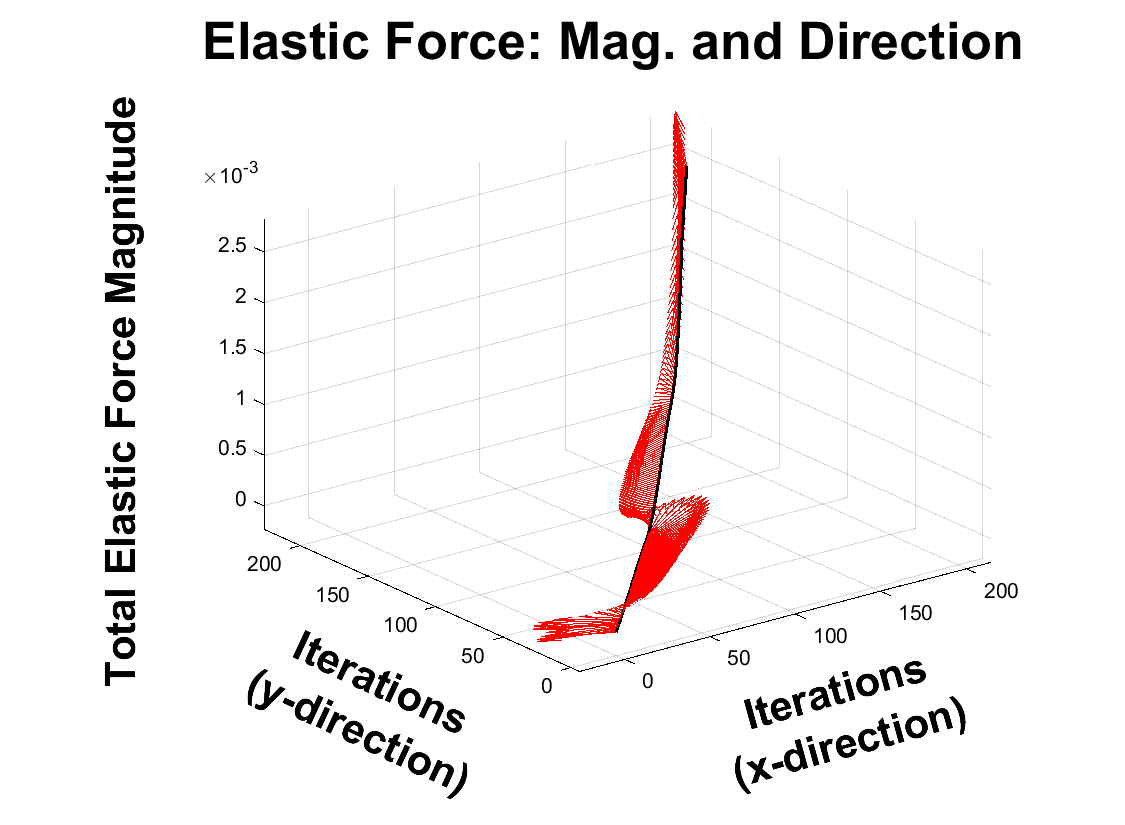}
\includegraphics[width=7cm]{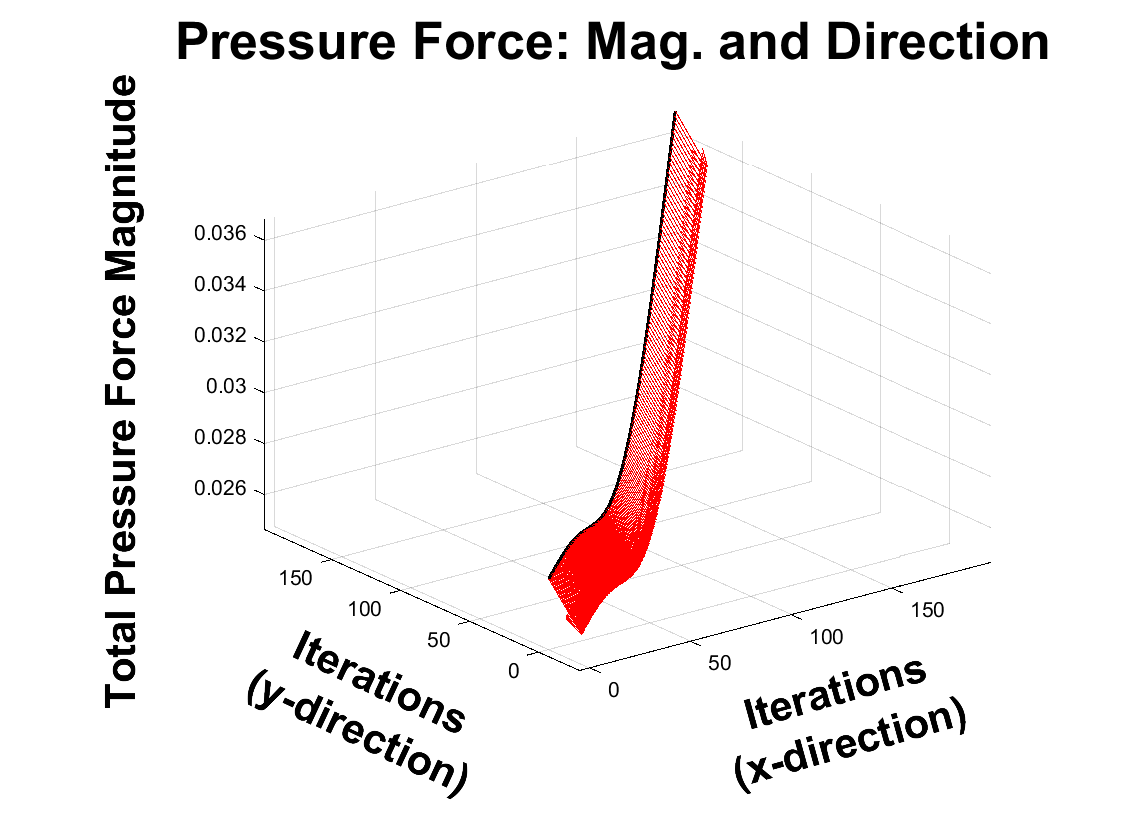}
\caption{$3D$ representations of total elastic (left) and pressure (right) force magnitudes and directions over iterations.}
\label{fig:ElasticPressureForce3D}
\end{center}
\end{figure}

To illustrate the evolution of force magnitudes over time, 3D visualizations of total pressure and elastic forces are provided in Figure \ref{fig:ElasticPressureForce3D}. 
These plots depict the force magnitudes along the simulation iterations (Z-axis) with arrows indicating the normalized force direction 
in $3D$ space at each time step (X- and Y-axes representing iterations).

\begin{figure}[H]
\begin{center}
\includegraphics[width=7cm]{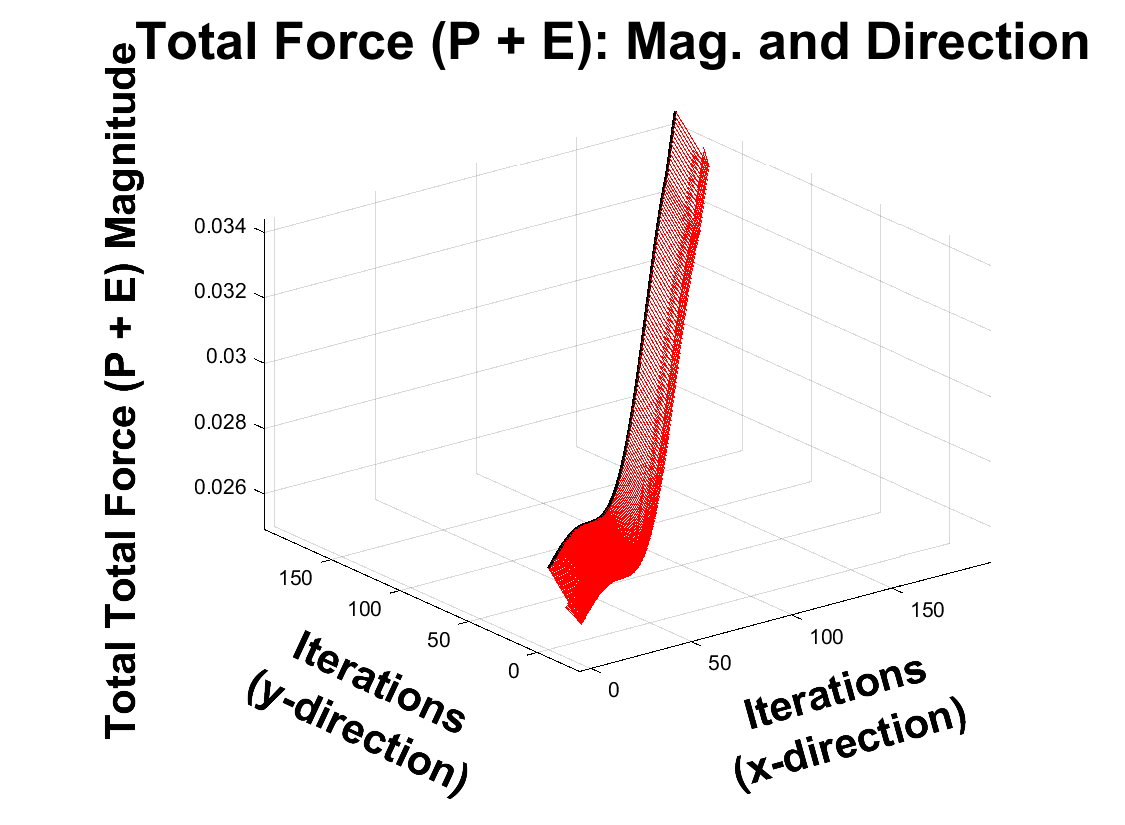}
\includegraphics[width=7cm]{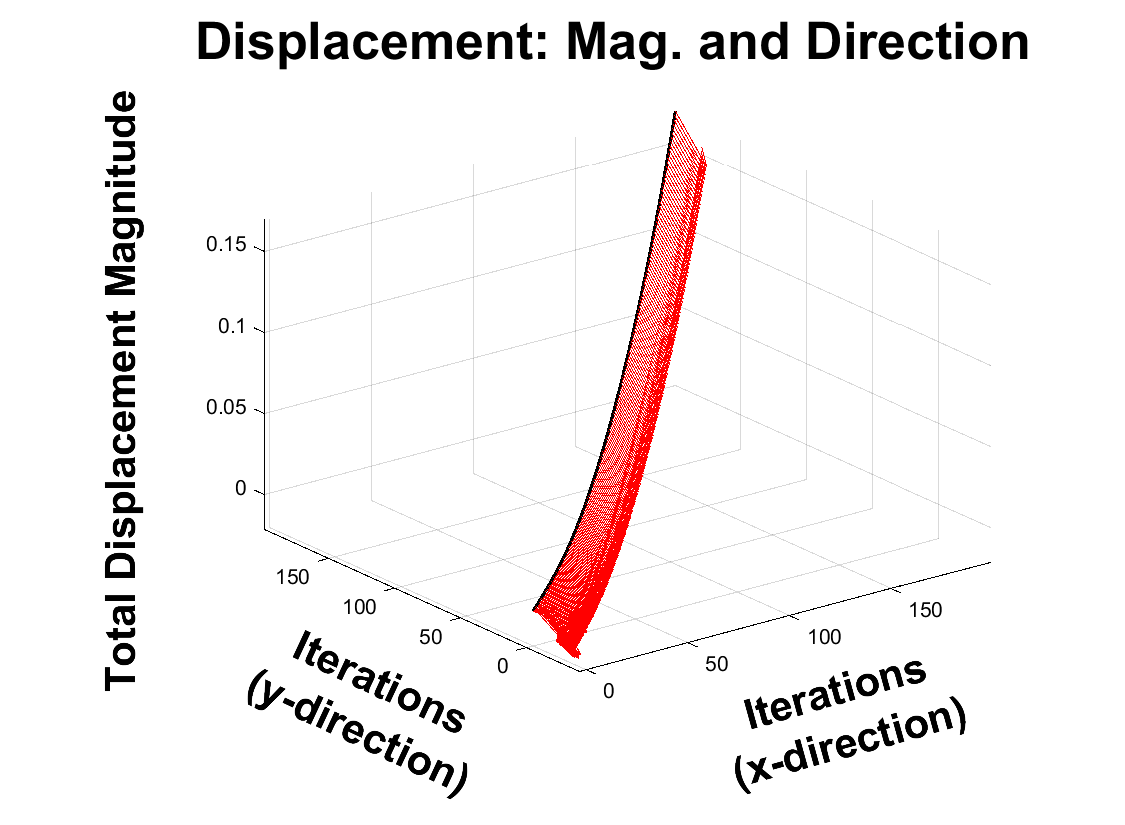}
\caption{Total force (elastic + pressure) exerted (left) and spatial displacement (right) of the test point over iterations.}
\label{fig:TotalForceDisplacement3D}
\end{center}
\end{figure}

Similarly, the total force (elastic + pressure) and the spatial displacement of the test point over iterations are 
illustrated in Figure \ref{fig:TotalForceDisplacement3D}.

\begin{figure}[H]
\begin{center}
\includegraphics[width=7cm]{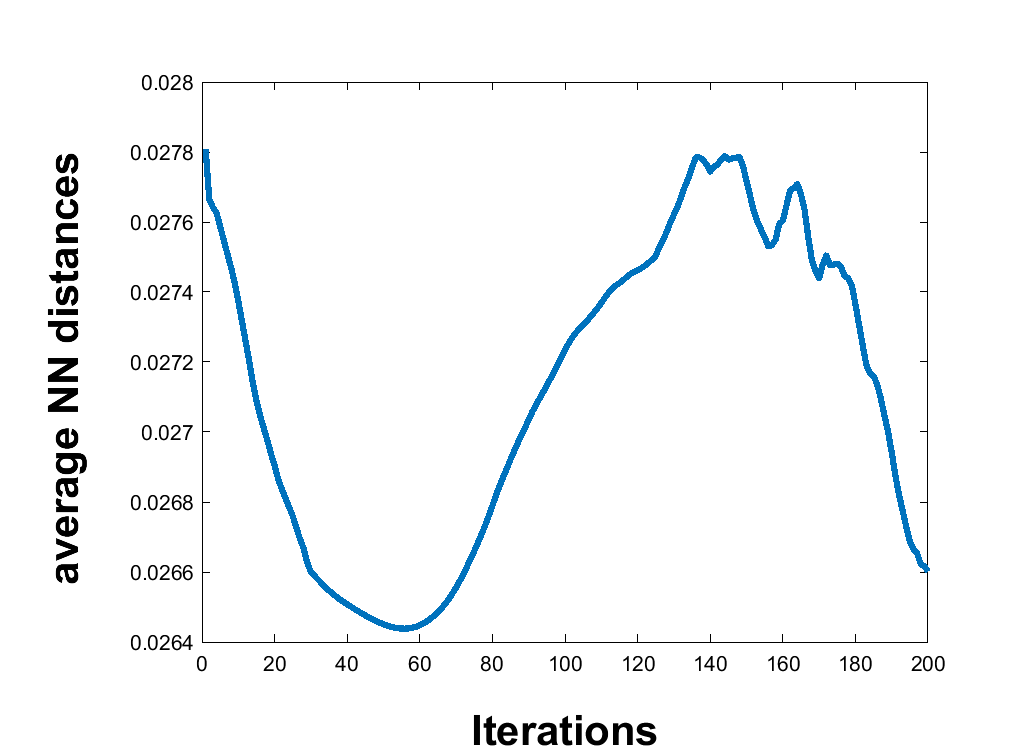}
\includegraphics[width=7cm]{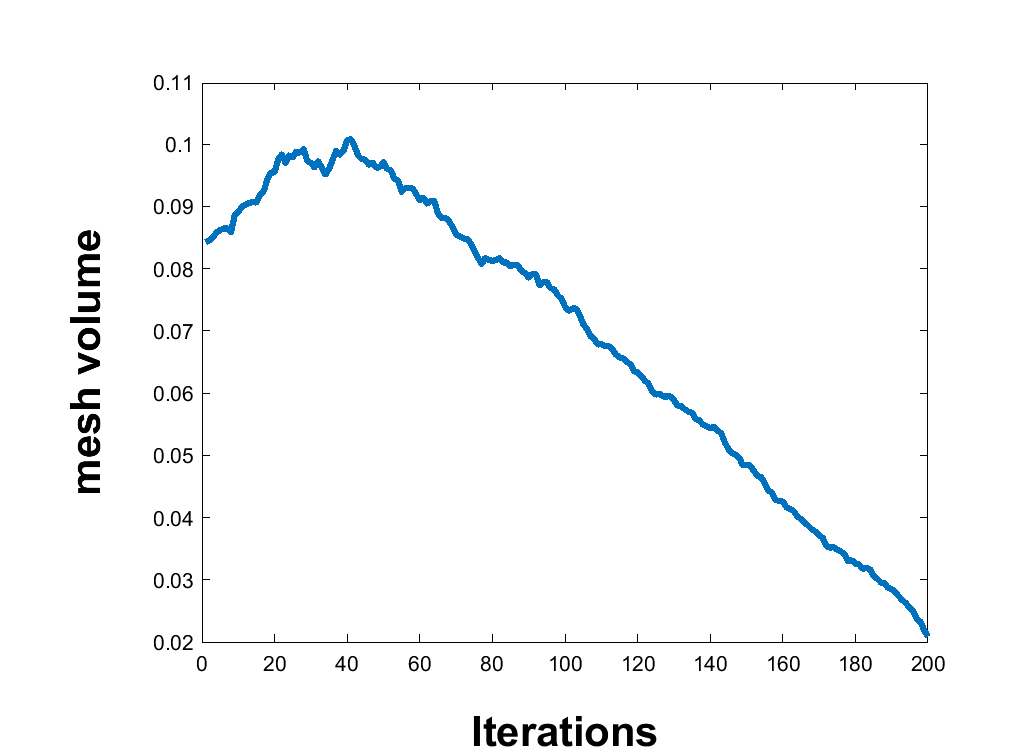}
\caption{Average nearest-neighbor distance (left) and inner volume of the mesh (right) as functions of simulation iterations.}
\label{fig:NN_volume}
\end{center}
\end{figure}

Additionally, the evolution of the average nearest-neighbor distance between adjacent mesh points and the change in inner volume of the closed mesh over the simulation iterations are shown in Figure \ref{fig:NN_volume}.

The results indicate that the deformation process successfully captures the concave characteristics of the focal body while preserving mesh 
closure and consistency. The iterative force dynamics, as depicted in Figures \ref{fig:PressureForceXYZ} and \ref{fig:ElasticForceXYZ}, 
show a progressive redistribution of elastic 
and pressure forces. Although these forces appear to stabilize temporarily, the system does not necessarily settle into a fully steady-state 
configuration. Instead, the evolution of the mesh continues, particularly in regions where local deformations persist due to ongoing force 
imbalances.

From a theoretical perspective, the simulation adheres to the principles of force balance between pressure-induced contraction and 
elasticity-driven restoration. Given that the pressure force dominates initially, the mesh undergoes substantial contraction before reaching 
a configuration where forces momentarily balance. However, due to the continuous adaptation of the mesh, the system does not reach a final 
steady-state but remains dynamically evolving. While this process reduces the total surface area, it does not strictly minimize surface 
area in the mathematical sense of a minimal surface.

An important observation is the relationship between mesh volume reduction and the nearest-neighbor distance evolution, as seen 
in Figure \ref{fig:NN_volume}. The contraction theorem discussed earlier (see also in \cite{Moriya:2024formconvexhullconcavity}) suggests that if pressure continuously exceeds elastic 
resistance, the mesh will eventually collapse to a near-zero volume state, with all cloud points lying on the surface. 
This aligns with the physical interpretation of an over-contracted flexible foil.

Furthermore, the $3D$ visualizations in Figures \ref{fig:ElasticPressureForce3D} and \ref{fig:TotalForceDisplacement3D} illustrate the directional consistency of force application over iterations. These results confirm that the method preserves force directionality, ensuring a controlled transition of mesh vertices towards their final configuration.

The findings validate the effectiveness of the FFMG approach for capturing concave structures. The behavior observed under different 
force conditions aligns with theoretical expectations, providing a strong basis for further refinement of the method, particularly 
in scenarios where the spatial distribution of the point cloud is characterized by enhanced concave nature.

\section{Conclusion}

This study introduces and applies the Flexible Foil Mesh Generation method to model the highly concave Focal Body
formed by the convergence of collimated light reflected from a spherical mirror. 
By leveraging a physically-driven mesh deformation process, we demonstrate that FFMG effectively reconstructs complex 
concave geometries that traditional computational geometry techniques struggle to capture. 

The theoretical foundation of FFMG rests on a force-equilibrium  damping model, where elastic forces preserve structural 
integrity, pressure-induced contraction drives surface evolution, and damping stabilizes the convergence process. 

Through iterative refinement, the simulated mesh dynamically adapts to the underlying structure of the FB, 
while ensuring enclosure of the spatially dense regions where reflected rays converge. 

This physically-motivated framework presents a significant advancement in surface reconstruction techniques, particularly for 
highly non-convex point clouds arising in optical, acoustical, and structural simulations.

The specific optical example analyzed, ray tracing of a spherical mirror, highlights the limitations of purely geometrical optics. 
While the present study focuses on ray intersection density as a proxy for light concentration, a truly comprehensive 
model of the FB would require incorporating wave optics principles such as diffraction and interference. 
Future refinements could integrate Huygens-Fresnel diffraction models, intensity-weighted ray tracing, or phase-space representations 
to better approximate the actual optical energy distribution within the FB. 

From a computational standpoint, our analysis of fixed vs. adaptive time-stepping demonstrates that the latter could significantly 
improve efficiency in the deformation process, especially in regions of high curvature variation. 
However, care must be taken to prevent overshooting and instabilities in regions where force distributions vary sharply.

The findings from this work are not restricted to optical simulations. The principles underlying FFMG offer a generalized methodology 
applicable to any domain requiring accurate reconstruction of concave surfaces, including medical imaging (e.g., organ 
boundary reconstruction), astrophysical simulations, fluid interface modeling, and inverse shape problems in engineering design. 
The ability to represent highly concave focal structures with minimal user intervention makes FFMG a promising tool for computational 
geometry and applied physics.

In conclusion, the effectiveness of FFMG in capturing concave spatial structures provides a foundation for future research, particularly 
in exploring force-driven surface evolution under varying physical constraints. 
The ability to intentionally manipulate equilibrium conditions in mesh contraction opens new avenues for optimizing geometrical 
representations of focal surfaces, wavefront interactions, and even dynamically evolving physical systems. 

Further development of hybrid geometric-waveform models could bridge the gap between purely geometric simulations and full optical-field 
representations, enhancing the accuracy and applicability of FFMG in high-precision domains.

\section*{Acknowledgments}
The author would like to express sincere gratitude to Dr. H. Primack for his significant contributions to 
the simulation work presented in this study. 
Dr. Primack's expertise in computational modeling and his assistance in developing the simulation algorithms 
were crucial to the successful execution of this research. 

\renewcommand{\bibname}{References}


\begin{thebibliography}{99}

\bibitem{2024Vaara}
Niklas Vaara, Pekka Sangi, Miguel Bordallo L\'opez, Janne Heikkil\"a,
\newblock A Ray Launching Approach for Computing Exact Paths with Point Clouds.
\newblock \emph{arXiv:2402.13747 [eess.SP]}, 2024.
\newblock DOI: \url{https://doi.org/10.48550/arXiv.2402.13747}.

\bibitem{klimes1994network}
Lud\v{e}k Klime\v{s} and Michal Kvasni\v{c}ka.
\newblock 3-D network ray tracing.
\newblock \emph{Geophysical Journal International}, 116(3):726--738, March 1994.
\newblock \url{https://doi.org/10.1111/j.1365-246X.1994.tb03293.x}.

\bibitem{Moriya:2025FlexableFoil}
Netzer Moriya,
\newblock Physically-Based Mesh Generation for Confined $3D$ Point Clouds Using Flexible Foil Models.
\newblock \emph{arXiv:2502.06541 [math.OC]}, 2025.
\newblock DOI: \url{https://doi.org/10.48550/arXiv.2502.06541}.

\bibitem{Moriya:2024TheLargest}
Netzer Moriya,
\newblock The Largest Empty Sphere Problem in $3D$ Hollowed Point Clouds.
\newblock \emph{arXiv:2401.07593 [math.OC]}, 2024.
\newblock DOI: \url{https://doi.org/10.48550/arXiv.2401.07593}.

\bibitem{Moriya:2024formconvexhullconcavity}
Netzer Moriya,
\newblock Form Convex Hull to Concavity: Surface Contraction Around a Point Set.
\newblock \emph{arXiv:2401.14189 [math.OC]}, 2024.
\newblock DOI: \url{https://doi.org/10.48550/arXiv.2401.14189}.

\bibitem[Loffler and van Kreveld(2010)]{Loffler:2010}
M. L\"offler and M. van Kreveld.
\newblock Largest and Smallest Convex Hulls for Imprecise Points.
\newblock \emph{Algorithmica}, 56(2):235--269, 2010.
\newblock \url{https://doi.org/10.1007/s00453-008-9174-2}
\newblock DOI: 10.1007/s00453-008-9174-2.

\bibitem{su1996delaunay}
P. Su and S.-M. Hu.
\newblock Three-dimensional Delaunay triangulation algorithms.
\newblock \emph{Advances in Engineering Software}, 27(1-2):55--60, 1996.

\bibitem{aurenhammer1991voronoi}
F. Aurenhammer.
\newblock Voronoi diagrams---a survey of a fundamental geometric data structure.
\newblock \emph{ACM Computing Surveys}, 23(3):345--405, 1991.

\bibitem{edelsbrunner1994alpha}
H. Edelsbrunner and E. P. M\"ucke.
\newblock Three-dimensional alpha shapes.
\newblock \emph{ACM Transactions on Graphics (TOG)}, 13(1):43--72, 1994.

\bibitem{fowles1975}
G. R. Fowles,
\newblock \emph{Introduction to Modern Optics}, 2nd ed., Dover Publications, 1975.

\bibitem{born1999} 
M. Born and E. Wolf, 
\newblock \emph{Principles of Optics}, 7th ed., Cambridge University Press, 1999.

\bibitem{berry1980}
M. V. Berry and C. Upstill, 
\newblock "Catastrophe Optics: Morphologies of Caustics and Their Diffraction Patterns" 
\newblock \emph{Progress in Optics}, vol. 18, pp. 257-346, 1980.

\bibitem{goodman2005}
J. W. Goodman, 
\newblock \emph{Introduction to Fourier Optics}, 3rd ed., McGraw-Hill, 2005.

\bibitem{arnold1992catastrophe}
V. I. Arnold.
\newblock \emph{Catastrophe Theory}.
\newblock Springer, 1992.

\bibitem{hecht2002} 
E. Hecht, 
\newblock \emph{Optics}, 4th ed., Addison-Wesley, 2002.

\bibitem{CLF1928} 
R. Courant, K. Friedrichs, and H. Lewy, 
\newblock \textit{\text{\"U}ber die partiellen Differenzengleichungen der mathematischen Physik}, Mathematische Annalen, 1928.

\end{thebibliography}
\end{document}